\documentclass[a4paper,11pt]{amsart}
\input{xypic}
\newtheorem{proposition}{Proposition}[section]
\newtheorem{lemma}[proposition]{Lemma}
\newtheorem{corollary}[proposition]{Corollary}
\newtheorem{theorem}[proposition]{Theorem}

\theoremstyle{definition}
\newtheorem{definition}[proposition]{Definition}
\newtheorem{example}[proposition]{Example}

\theoremstyle{remark}
\newtheorem{remark}[proposition]{Remark}

\newcommand{\thlabel}[1]{\label{th:#1}}
\newcommand{\thref}[1]{Theorem~\ref{th:#1}}
\newcommand{\selabel}[1]{\label{se:#1}}
\newcommand{\seref}[1]{Section~\ref{se:#1}}
\newcommand{\lelabel}[1]{\label{le:#1}}
\newcommand{\leref}[1]{Lemma~\ref{le:#1}}
\newcommand{\prlabel}[1]{\label{pr:#1}}
\newcommand{\prref}[1]{Proposition~\ref{pr:#1}}
\newcommand{\colabel}[1]{\label{co:#1}}
\newcommand{\coref}[1]{Corollary~\ref{co:#1}}

\newcommand{\relabel}[1]{\label{re:#1}}
\newcommand{\reref}[1]{Remark~\ref{re:#1}}
\newcommand{\exlabel}[1]{\label{ex:#1}}
\newcommand{\exref}[1]{Example~\ref{ex:#1}}
\newcommand{\delabel}[1]{\label{de:#1}}

\newcommand{\eqlabel}[1]{\label{eq:#1}}
\newcommand{\equref}[1]{(\ref{eq:#1})}

\def\equal#1{\smash{\mathop{=}\limits^{#1}}}

\newcommand{\can}{{\rm can}}

\newcommand{\Hom}{{\rm Hom}}
\newcommand{\HOM}{{\rm HOM}}
\newcommand{\End}{{\rm End}}
\newcommand{\END}{{\rm END}}

\newcommand{\Ker}{{\rm Ker}\,}

\def\lan{\langle}
\def\ran{\rangle}
\def\ot{\otimes}

\def\sq{\square}

\def\rightactr{\hbox{$\leftharpoonup$}}

\newcommand{\Cc}{\mathcal{C}}
\newcommand{\Dd}{\mathcal{D}}
\newcommand{\Ee}{\mathcal{E}}

\newcommand{\Mm}{\mathcal{M}}

\def\*C{{}^*\hspace*{-1pt}{\Cc}}
\def\text#1{{\rm {\rm #1}}}

\def\Morita{\dul{\rm Morita}}

\def\ol{\overline}
\def\ul{\underline}
\def\dul#1{\underline{\underline{#1}}}

\def\mapright#1{\smash{\mathop{\longrightarrow}\limits^{#1}}}

\def\sf{\scriptstyle}

\begin{document}

\title[Bimodules over Hopf-Galois extensions]{Morita
equivalences induced by bimodules over Hopf-Galois extensions}
\author[S. Caenepeel]{Stefaan Caenepeel}
\address{Faculty of Engineering,
Vrije Universiteit Brussel, VUB, Pleinlaan 2, B-1050 Brussels, Belgium}
\email{scaenepe@vub.ac.be}
\urladdr{http://homepages.vub.ac.be/\~{}scaenepe/}

\author[S. Crivei]{Septimiu Crivei}
\address{Faculty of Mathematics and Computer Science,
Babe\c s-Bolyai University, Str. Mihail Kog\u alniceanu 1,
RO-400084 Cluj-Napoca, Romania} \email{crivei@math.ubbcluj.ro}

\author[A. Marcus]{Andrei Marcus}
\address{Faculty of Mathematics and Computer Science,
Babe\c s-Bolyai University, Str. Mihail Kog\u alniceanu 1,
RO-400084 Cluj-Napoca, Romania} \email{marcus@math.ubbcluj.ro}

\author[M. Takeuchi]{Mitsuhiro Takeuchi}
\address{Institute of Mathematics, University of Tsukuba, Tsukuba, Ibaraki 305, Japan}
\email{takeu@math.tsukuba.ac.jp}


\begin{abstract}
Let $H$ be a Hopf algebra, and $A,B$ be $H$-Galois
extensions. We investigate the category ${}_A\mathcal{M}_B^H$ of relative Hopf bimodules,
and the Morita equivalences between $A$ and $B$ induced by them.
\end{abstract}

\thanks{This research was supported by the bilateral project BWS04/04 ``New
Techniques in Hopf algebras and graded ring theory" of the Flemish and Romanian
governments and by the research project G.0622.06 ``Deformation quantization methods
for algebras and categories with applications to quantum mechanics" from
FWO-Vlaanderen.
The second author  acknowledges the support of the grant CEEX-ET 47/2006.
The third author acknowledges the support of a Bolyai Fellowship
of the Hungarian Academy of Science}

\maketitle

\section*{Introduction}
This paper is a contribution to the representation theory of
Hopf-Galois extensions, as originated by Schneider in
\cite{Schneider1}. More specifically, we consider the following
questions. Let $H$ be a Hopf algebra, and $A$, $B$ right
$H$-comodule algebras. Moreover, assume that $A$ and $B$ are
right faithfully flat $H$-Galois extensions.
\begin{enumerate}
\item If $A$ and $B$ are Morita equivalent, does it follow that
$A^{{\rm co}H}$ and $B^{{\rm co}H}$ are also Morita equivalent?
\item Conversely,
if $A^{{\rm co}H}$ and $B^{{\rm co}H}$ are Morita equivalent, when does it follow
that $A$ and $B$ are Morita equivalent?
\end{enumerate}
These questions have been considered in \cite{Marcus} in the context of
strongly group graded algebras, the motivation coming from problems raised in
the modular representation theory of finite groups. The results of the present
paper generalize the results of \cite[Sections 2 and 3]{Marcus}.\\
Given a right $H$-comodule algebra $A$, and a left $H$-comodule
algebra $B$, we consider $(A\ot B,H)$-Hopf modules. These are at
the same time left $A\ot B$-modules and right $H$-comodules, with
a suitable compatibility condition. There are various ways to look
at these Hopf modules: they are Doi-Hopf modules (see \cite{Doi})
over a certain Doi-Hopf datum (with two possible descriptions of
the underlying module coalgebra), and they can also be viewed as
comodules over a coring (see \seref{3}). The main result of
\seref{2}, and also the main tool used during the rest of the
paper, is a structure Theorem for $(A\ot B,H)$-Hopf modules,
stating that the category of $(A\ot B,H)$-Hopf modules is
equivalent to the category of left modules over the cotensor
product $A\square_H B$, under the condition that $A$ is a
faithfully
flat $H$-Galois extension.\\
The results from \seref{2} can be applied to relative Hopf bimodules: let $A$ and
$B$ be right $H$-comodule algebras, and consider $(A,B)$-bimodules with a
right $H$-coaction, satisfying a certain compatibility condition. The category
of relative Hopf bimodules is then isomorphic to the category of
$(A\ot B^{\rm op},H)$-Hopf modules. In \seref{4}, we state the Structure Theorem
for relative Hopf bimodules, and  we investigate the compatibility of
the category equivalence with the Hom and tensor functors.\\
In \seref{5}, we apply our results to discuss the two problems
stated above. We introduce the notion of $H$-Morita contexts, and
we show that if two right faithfully flat $H$-Galois extensions are
connected by a (strict) $H$-Morita context, then the algebras of
coinvariants are also connected by a (strict) Morita context. Our
main result is the following converse result: if the algebras of
coinvariants are Morita equivalent, in such a way that the bimodule
structure on one of the connecting modules can be extended to a
left-action by the cotensor product $A\square_H B^{\rm op}$, then
$A$ and $B$ are $H$-Morita
equivalent.\\
In \seref{6}, we show that the Morita equivalence coming from a strict $H$-Morita
context between two faithfully flat $H$-Galois extensions respects the Miyashita-Ulbrich
action.
In \seref{7}, we investigate the behavior of $H$-Morita equivalences with respect to Hopf subalgebras. \\
The category of relative Hopf modules and $A$-linear (not necessarily $H$-colinear)
modules is an $H$-colinear category. If two right $H$-comodule algebras are
$H$-Morita equivalent, then the induced equivalence between their categories of
relative Hopf modules is $H$-colinear. In \seref{8}, we study the converse property:
when does every $H$-colinear equivalence between two categories of relative Hopf
modules come from a strict $H$-Morita context. This leads to a generalization of
the Eilenberg-Watts Theorem (\prref{8.5}). The main result is \coref{8.7}, stating that
every $H$-colinear equivalence comes from a strict $H$-Morita context if the 
Hopf algebra $H$ is projective, and the $H$-comodule algebras
$A$ and $B$ are $H$-Galois extensions of their subalgebras of coinvariants.\\
For basic results on Hopf algebras, we refer the reader to
\cite{DNR} or \cite{Mont}. For a concise treatment of corings and their applications,
we refer to \cite{BrzezinskiWisbauer}.

\section{Preliminary results}\selabel{1}
Throughout this paper $H$ is a Hopf algebra over a commutative
ring $k$, with bijective antipode $S$. We use the Sweedler
notation for the comultiplication on $H$: $\Delta(h)= h_{(1)}\ot
h_{(2)}$. $\Mm^H$ (respectively ${}^H\Mm$) is the category of
right (respectively left) $H$-comodules. For a right $H$-coaction
$\rho$ (respectively a left $H$-coaction $\lambda$) on a
$k$-module $M$, we denote
$$\rho(m)=m_{[0]}\ot m_{[1]}\quad{\rm and}\quad\lambda(m)=m_{[-1]}\ot m_{[0]}.$$
The submodule of coinvariants $M^{{\rm co}H}$ of a right
(respectively left) $H$-co\-mo\-du\-le $M$ consists of the elements
$m\in M$ satisfying $\rho(m)=m\ot 1$ (respectively
$\lambda(m)=1\ot m$).\\
Let $A$ be a right $H$-comodule algebra. ${}_A\Mm^H$ and $\Mm^H_A$ are
the categories of left and right relative Hopf modules. We have two pairs of
adjoint functors $(F_1=A\ot_{A^{{\rm co}H}}-,\ G_1=(-)^{{\rm co}H})$
and $(F_2=-\ot_{A^{{\rm co}H}}A,\ G_2=(-)^{{\rm co}H})$ between the categories
${}_{A^{{\rm co}H}}\Mm$ and ${}_A\Mm^H$, and between
$\Mm_{A^{{\rm co}H}}$ and $\Mm^H_A$. The unit and counit of the adjunction
$(F_1,G_1)$ are given by the formulas
$$\eta_{1,N}:\ N\to (A\ot_{A^{{\rm co}H}} N)^{{\rm co}H},\quad\eta_{1,N}(n)=1\ot n;$$
$$\varepsilon_{1,M}:\ A\ot_{A^{{\rm co}H}} M^{{\rm co}H}\to M,\quad\varepsilon_{1,M}(a\ot m)=am.$$
The formulas for the unit and counit of $(F_2,G_2)$ are similar. Consider the canonical maps
$$\can:\ A\ot_{A^{{\rm co}H}} A\to A\ot H,~~\can(a\ot b)=ab_{[0]}\ot b_{[1]};$$
$$\can':\ A\ot_{A^{{\rm co}H}} A\to A\ot H,~~\can(a\ot b)=a_{[0]}b\ot a_{[1]}.$$
It is well-known (see for example \cite{KT}) that $\can$ is an isomorphism if and only if
$\can'$ is an isomorphism.

\begin{theorem}\thlabel{1.1}Let $A$ be a right $H$-comodule algebra. 
Consider the following statements:
\begin{enumerate}
\item $(F_2,G_2)$ is a pair of inverse equivalences;
\item $(F_2,G_2)$ is a pair of inverse equivalences and $A\in {}_{A^{{\rm co}H}}\Mm$ is flat;
\item $\can$ is an isomorphism and $A\in {}_{A^{{\rm co}H}}\Mm$ is faithfully flat;
\item $(F_1,G_1)$ is a pair of inverse equivalences;
\item $(F_1,G_1)$ is a pair of inverse equivalences and $A\in \Mm_{A^{{\rm co}H}}$ is flat;
\item $\can'$ is an isomorphism and $A\in \Mm_{A^{{\rm co}H}}$ is faithfully flat.
\end{enumerate}
We have the following implications:
$$(3)\Longleftrightarrow (2)\Longrightarrow (1)~~;~~(6)\Longleftrightarrow (5)\Longrightarrow (4).$$
If $H$ is flat as a $k$-module, then $(1)\Longleftrightarrow (2)$ and $(4)\Longleftrightarrow (5)$.\\
If $k$ is field then the six conditions are equivalent.
\end{theorem}

If the first three conditions of \thref{1.1} hold, then we call $A$
a left faithfully flat $H$-Galois extension; if the three other conditions hold, then we call
$A$ a right faithfully flat $H$-Galois extension.

\begin{proof}
The equivalence of (2) and (3) is well-known. It is essentially \cite[Theorem 3.7]{Schneider},
which is an improvement of \cite[Theorem 2.11]{DT}.
For the equivalence of (5) and (6), we observe that $A$ is a left $H^{\rm cop}$-comodule,
so, by the left handed version of the equivalence $(3)\Longleftrightarrow (2)$,
(6) is equivalent to flatness of $A\in \Mm_{A^{{\rm co}H}}$ and equivalence between
the categories ${}_{A^{{\rm co}H}}\Mm$ and ${}_A^{H^{\rm cop}}\Mm\cong {}_A\Mm^H$.\\
The implications $(2)\Longrightarrow (1)$ and $(5)\Longrightarrow (4)$ are trivial.\\
If $H$ is flat as a $k$-module, then $\Mm_A^H$ is an abelian category and the
forgetful functor $\Mm_A^H\to \Mm_A$ is exact. 
If $F_2$ is an equivalence, then the functor $-\ot_{A^{{\rm co}H}}:\
\Mm_{A^{{\rm co}H}}\to \Mm_A$ is exact since it is the composition of the forgetful
functor and the equivalence $F_2$. This shows that $A$ is flat as a left
$A^{{\rm co}H}$-module, and the implication $(1)\Longrightarrow (2)$ follows.
$(4)\Longrightarrow (5)$ can be proved in a similar way.\\
If $k$ is a field, then the equivalence of the six statements in the Theorem follows
from \cite[Theorem I]{Schneider0}.
\end{proof}

Let $M$ be a right $H$-comodule, and $N$ a left $H$-comodule. The cotensor
product $M\square_H N$ is the $k$-module
$$M\square_H N= \{\sum_i m_i\ot n_i\in M\ot N~|~\sum_i \rho(m_i)\ot n_i=
\sum_i m_i\ot \lambda(n_i)\}.$$
If $H$ is cocommutative, then $M\square_H N$ is also a right (or left) $H$-comodule.

\begin{proposition}\prlabel{1.2}
Let $R$ be a $k$-algebra, and assume that $P\in \Mm_R$ is flat.
Take $M\in {}_R\Mm^H$
and $N\in {}^H\Mm$, and assume that we have a right $H$-coaction on $M$ that
is left $R$-linear. Then the map
$$P\ot_R (M\square_H N)\to (P\ot_R M)\square_H N,~~
p\ot (\sum_i m_i\ot n_i)\mapsto \sum_i (p\ot m_i)\ot n_i$$
is bijective.
\end{proposition}

\section{A structure theorem for $(A\ot B,H)$-Hopf modules}\selabel{2}
Under our assumption on $H$, $H\ot H^{\rm cop}$ is also a Hopf
algebra, and $H$ is a left $H\ot H^{\rm cop}$-module coalgebra;
the left $H\ot H^{\rm cop}$-action is given by
$$(k\ot l)\cdot h=khS(l),$$
for all $h,k,l\in H$.\\

We present an alternative description of $H$ as a left $H\ot H^{\rm cop}$-module coalgebra.
$H\ot H^{\rm cop}\in {}_{H\ot H^{\rm cop}}\Mm_H$, with right
$H$-action induced by the comultiplication on $H$, and $k\in {}_H\Mm$ via $\varepsilon$,
so we have the left $H\ot H^{\rm cop}$-module $(H\ot H^{\rm cop})\ot_H k$.
$(H\ot H^{\rm cop})\ot_H k$ is a coalgebra with comultiplication and counit given by
$$\Delta((h\ot h')\ot_H 1)=(h_{(1)}\ot h'_{(2)})\ot_H 1\ot (h_{(2)}\ot h'_{(1)})\ot_H 1;$$
$$\varepsilon((h\ot h')\ot_H 1)=\varepsilon(hh').$$
It is easy to show that $(H\ot H^{\rm cop})\ot_H k$ is an $H\ot H^{\rm cop}$-module coalgebra.

\begin{proposition}\prlabel{2.0}
$(H\ot H^{\rm cop})\ot_H k$ and $H$ are isomorphic as
$H\ot H^{\rm cop}$-module coalgebras.
\end{proposition}

\begin{proof}
Define
$$f:\ (H\ot H^{\rm cop})\ot_H k\to H,~~f((h\ot h')\ot_H 1)=hS(h');$$
$$g:\ H\to (H\ot H^{\rm cop})\ot_H k,~~g(h)=(h\ot 1)\ot_H 1.$$
$f$ is well-defined since for all $h,h',l\in H$
$$f((h\ot h')l\ot_H 1)=hl_{(1)}S(h'l_{(2)})=hS(h')\varepsilon(l)=
f((h\ot h')\ot_H \varepsilon(l)).$$
$f$ is $H\ot H^{\rm cop}$-linear since  for all $h,h',k,k'\in H$
$$f((kh\ot k'h')\ot_H 1)=khS(k'h')=(k\ot k')\cdot (hS(k'))=(k\ot k')f((h\ot h')\ot_H 1).$$
$f$ is a coalgebra map since for all $h,h'\in H$
$$((f\ot f)\circ\Delta)((h\ot h')\ot_H 1)=h_{(1)}S(h'_{(2)})\ot h_{(2)}S(h'_{(1)})=
\Delta(hS(h')),$$
and
$$(\varepsilon\circ f)((h\ot h')\ot_H 1)=\varepsilon(hS(h'))=\varepsilon(hh').$$
It is obvious that $f\circ g=H$. Finally  for all $h,k\in H$
\begin{eqnarray*}
&&\hspace*{-2cm}
(g\circ f)((h\ot k)\ot_H 1)= g(hS(k))\\
&=&(hS(k)\ot 1)\ot_H 1
= (hS(k_{(1)})\ot 1)\ot_H\varepsilon(k_{(2)})\\
&=&
(hS(k_{(1)})k_{(2)}\ot k_{(3)})\ot_H 1=(h\ot k)\ot_H 1.
\end{eqnarray*}
\end{proof}

Let $A$ be a right $H$-comodule algebra, and $B$ a left
$H$-comodule algebra. Then $A\ot B$ is a right $H\ot H^{\rm
cop}$-comodule algebra, with coaction
$$\rho(a\ot b)=a_{[0]}\ot b_{[0]}\ot a_{[1]}\ot b_{[-1]}.$$
Then $(H\ot H^{\rm cop}, A\ot B,H)$ is a left-right Doi-Hopf datum (see
\cite{CMZ} or \cite{Doi} for details), and we can consider the category
${}_{A\ot B}\Mm(H\ot H^{\rm cop})^H$ of Doi-Hopf modules. The objects of this
category are $k$-modules $M$ with a left $A\ot B$-action and a right $H$-coaction such
that
$$\rho((a\ot b)m)=(a_{[0]}\ot b_{[0]})m_{[0]} \ot a_{[1]}m_{[1]}S(b_{[-1]}),$$
for all $a\in A$, $b\in B$ and $m\in M$.
The objects of ${}_{A\ot B}\Mm(H\ot H^{\rm cop})^H$ are called $(A\ot B,H)$-Hopf modules.
It is well-known and easily verified that  $A\ot B\in {}_{A\ot B}\Mm(H\ot H^{\rm cop})^H$,
with coaction defined by
$$\rho(a\ot b)= a_{[0]}\ot b_{[0]}\ot a_{[1]}S(b_{[-1]}).$$

\begin{lemma}\lelabel{2.1}
With notation as above, we have that $(A\ot B)^{{\rm co}H}=A\square_H B$.
\end{lemma}

\begin{proof}
Take $x=\sum_i a_i\ot b_i \in (A\ot B)^{{\rm co}H}$. Then
$$\sum_i a_i\ot b_i\ot 1= \sum_i a_{i[0]}\ot b_{i[0]}\ot a_{i[1]}S(b_{i[-1]}).$$
Apply $\lambda$ to the second tensor factor. Then switch the second and fourth
tensor factor, and multiply the third and fourth tensor factor. It follows that
$$\sum_i a_i\ot b_{i[0]}\ot b_{i[-1]}=\sum_i a_{i[0]}\ot b_{i[0]}\ot a_{i[1]}S(b_{i[-2]})b_{i[-1]}
=\sum_i a_{i[0]}\ot b_i \ot a_{i[1]},$$ and then $x\in A\square_H
B$. The converse inclusion is proved in a similar way.
\end{proof}

Recall (see for example \cite{Doi}) that we have a pair of adjoint functors $(F,G)$:
$$F:\ {}_{A\square_H B}\Mm\to {}_{A\ot B}\Mm(H\ot H^{\rm cop})^H,~~
F(N)=(A\ot B)\ot_{A\square_H B} N;$$
$$G:\  {}_{A\ot B}\Mm(H\ot H^{\rm cop})^H\to {}_{A\square_H B}\Mm,~~
G(M)=M^{{\rm co}H}.$$
The unit and counit of the adjunction are the following:
$$\eta_N:\ N\to \bigl((A\ot B)\ot_{A\square_H B} N\bigr)^{{\rm co}H},~~
\eta_N(n)=1_A\ot 1_B\ot n;$$
$$\varepsilon_M:\ (A\ot B)\ot_{A\square_H B} M^{{\rm co}H}\to M,~~
\varepsilon_M(a\ot b\ot m)=(a\ot b)m.$$

\begin{proposition}\prlabel{2.2}
Assume that $H$ is flat as a $k$-algebra. Let $A$ be a right
$H$-comodule algebra, and $B$ a left $H$-comodule algebra. We have
a right $H$-colinear map
$$f:\ A\ot_{A^{{\rm co}H}}(A\square_HB)=F_1(A\square_H B)\to A\ot B,~~
f(a\ot(\sum_i a_i\ot b_i))=\sum_i aa_i\ot b_i.$$
If $A$ is a right faithfully flat $H$-Galois extension, then $f$ is an isomorphism.
\end{proposition}

\begin{proof}
$f$ is right $H$-colinear since
\begin{eqnarray*}
&&\hspace*{-2cm}
\rho\bigl(f(a\ot(\sum_i a_i\ot b_i))\bigr)=
\sum_i a_{[0]}a_{i[0]}\ot b_{i[0]}\ot a_{[1]}a_{i[1]}S(b_{i[-1]})\\
&=& \sum_i a_{[0]}a_{i}\ot b_{i[0]}\ot a_{[1]}b_{i[-2]}S(b_{i[-1]})\\
&=&\sum_i a_{[0]}a_{i}\ot b_i\ot a_{[1]}\\
&=& (f\ot H)(\rho(a\ot(\sum_i a_i\ot b_i))).
\end{eqnarray*}
On $A\ot_{A^{{\rm co}H}}A$ and $A\ot H$, we consider the following right $H$-coactions:
$$\rho(a\ot b)=a\ot b_{[0]}\ot b_{[1]}~~;~~\rho(a\ot h)=a\ot h_{(1)}\ot h_{(2)}.$$
Then $\can:\ A\ot_{A^{{\rm co}H}}A\to A\ot H$ is right $H$-colinear, so we can
consider the map $\can\square_H B: (A\ot_{A^{{\rm co}H}}A)\square_H B\to (A\ot H)
\square_H B$. If $A$ is a right faithfully flat $H$-Galois extension, then $\can\square_H B$
is bijective, and applying \prref{1.2}, we see that $f$ is the composition of the
following isomorphisms:
\begin{eqnarray*}
&&\hspace*{-2cm}
A\ot_{A^{{\rm co}H}}(A\square_HB)\cong (A\ot_{A^{{\rm co}H}}A)\square_HB\\
&\cong&(A\ot H)\square_H B\cong A\ot (H\square_H B)\cong A\ot B.
\end{eqnarray*}
\end{proof}

The following structure theorem is the main result of this
section.

\begin{theorem}\thlabel{2.3} Let $A$ be a right $H$-comodule algebra, and $B$ a left
$H$-comodule algebra. If $A$ is a right faithfully flat $H$-Galois
extension, then $(F,G)$ is a pair of inverse equivalences between
the categories ${}_{A\square_H B}\Mm$ and ${}_{A\ot B}\Mm(H\ot
H^{\rm cop})^H$.
\end{theorem}

\begin{proof}
Take $N\in {}_{A\square_H B}\Mm$. We have a well-defined algebra map
$A^{{\rm co}H}\to A\square_H B$, sending $a$ to $a\ot 1_B$, and $N$ is a
left $A^{{\rm co}H}$-module, by restriction of scalars. Consider the isomorphism
\begin{eqnarray*}
&&\hspace*{-1cm}
\alpha_N=f\ot_{A\square_H B} N:\
F_1(N)=A\ot_{A^{{\rm co}H}} N\cong A\ot_{A^{{\rm co}H}} (A\square_H B)
\ot_{A\square_H B} N\\
&\to& F(N)= (A\ot B)\ot_{A\square_H B} N.
\end{eqnarray*}
It is easy to see that $\alpha_N(a\ot n)=(a\ot 1)\ot_{A\square_H B} n$, and $\alpha_N$
is right $H$-colinear since
$$(\alpha_N\ot H)((a_{[0]}\ot n)\ot a_{[1]})=((a_{[0]}\ot 1)\ot_{A\square_H B}n)\ot a_{[1]}
=\rho((a\ot 1)\ot n).$$
It follows that $\alpha_N$ restricts to  an isomorphism
$$\alpha_N^{{\rm co}H}:\ (A\ot_{A^{{\rm co}H}} N)^{{\rm co}H}\to
((A\ot B)\ot_{A\square_H B} N)^{{\rm co}H}.$$
It is then easily seen that
$$\eta_N=\alpha_N^{{\rm co}H}\circ \eta_{1,N}.$$
$\eta_{1,N}$ is an isomorphism by \thref{1.1}, and it follows that $\eta_N$ is
an isomorphism.\\
Take $M\in {}_{A\ot B}\Mm(H\ot H^{\rm cop})^H$. Then $M$ is a left $A$-module,
by restriction of scalars, and a relative Hopf module since
$$\rho(am)=\rho((a\ot 1)m)=(a_{[0]}\ot 1)m_{[0]}\ot a_{[1]}m_{[1]}S(1)=
a_{[0]}m_{[0]}\ot a_{[1]}m_{[1]}.$$
It is then easy to see that
$$\varepsilon_M\circ \alpha_{M^{{\rm co}H}}=\varepsilon_{1,M}.$$
It follows from \thref{1.1} that $\varepsilon_{1,M}$ is an isomorphism, and this
implies that $\varepsilon_M$ is an isomorphism.
\end{proof}

\section{Connection to comodules over corings}\selabel{3}
Let $A$ be a ring. Recall that an $A$-coring $\Cc$ is a comonoid in the monoidal
category ${}_A\Mm_A$. For a detailed discussion of the theory of corings and
comodules, we refer to \cite{BrzezinskiWisbauer}. One of the results is that we
can associate a coring to a Doi-Hopf datum, and that the category of Doi-Hopf modules
is isomorphic to the category of comodules over this coring.\\
Let us describe the $A\ot B$-coring $\Cc$ associated to the
left-right Doi-Hopf datum $(H\ot H^{\rm cop}, A\ot B, H)$ that we
have discussed in the previous section. We have that $\Cc=H\ot
A\ot B$, with left and right $A\ot B$-action given by
$$(a'\ot b')(h\ot a\ot b)(a''\ot b'')=
a'_{[1]}hS(b'_{[-1]})\ot a'_{[0]}aa''\ot b'_{[0]}bb''.$$
The comultiplication and counit are given by the formulas
$$\Delta(h\ot a\ot b)=(h_{(2)}\ot 1_A\ot 1_B)\ot_{A\ot B} (h_{(1)}\ot a\ot b);$$
$$\varepsilon(h\ot a\ot b)=\varepsilon(h)a\ot b.$$
The category ${}^\Cc\Mm$ of left $\Cc$-comodules is isomorphic to
${}_{A\ot B}\Mm(H\ot H^{\rm cop})^H$.\\
A Galois theory for corings can be developed (see
\cite{Brzezinski02,Caenepeel03}). Let $x$ be a grouplike element
of a coring $\Cc$, and let $$A^{{\rm co}\Cc}= \{a\in
A~|~ax=xa\}.$$ Then we have an adjoint pair of functors between
${}_{A^{{\rm co}\Cc}}\Mm$ and ${}^\Cc\Mm$. If this adjoint pair is
a pair of inverse equivalences, then the map
$$\can:\ A\ot_{A^{{\rm co}\Cc}}A\to \Cc,\quad\can(a\ot b)=axb$$
is an isomorphism of corings (see \cite[Proposition
3.1]{Caenepeel03}). If, in addition, $A$ is flat as a right
$A^{{\rm co}\Cc}$-module, then it also follows that $A$ is
faithfully flat as a right $A^{{\rm co}\Cc}$-module (see
\cite[Proposition 3.8, $2)\Rightarrow 1)$]{Caenepeel03}). We will apply
this to the coring $\Cc=H\ot A\ot B$. $1_H\ot 1_A\ot 1_B$ is a grouplike element of $H\ot A\ot B$, and
the associated pair of adjoint functors is precisely $(F,G)$. It
can be easily verified that the corresponding canonical map is
precisely the map
\begin{eqnarray}
&&\hspace*{-2cm}
\can:\ (A\ot B)\ot_{A\square_H B}(A\ot B)\to H\ot A\ot B ,\nonumber\\
&&\can((a\ot b)\ot (a'\ot b'))=a_{[1]}S(b_{[-1]})\ot a_{[0]}a'\ot b_{[0]}b'\eqlabel{coring1}
\end{eqnarray}

\begin{proposition}\prlabel{coring}
Let $A$ be a right $H$-comodule algebra, and $B$ a left
$H$-comodule algebra. Assume that $A$ is a right faithfully flat
$H$-Galois extension. Then $\can$
is an isomorphism. Furthermore, $A\ot B$ is faithfully flat as a right
$A\square_H B$-module.
\end{proposition}

\begin{proof}
It follows from \thref{2.3} that $(F,G)$ is a pair of inverse equivalences, hence
$\can$ is an isomorphism.\\
We will now show that $A\ot B$ is flat as a right
$A\square_H B$-module.
Assume that $N\to N'$ is a monomorphism of left $A\square_H B$-modules.
Using \prref{2.2} and the fact that $A$ is flat as a right $A^{{\rm co}H}$-module,
we find that
$$(A\ot B)_{A\square_HB} N\cong A\ot _{A^{{\rm co}H}} N
\to A\ot _{A^{{\rm co}H}}N'\cong (A\ot B)_{A\square_HB} N'$$
is injective. As explained above, it then follows from  \cite[Proposition 3.8]{Caenepeel03} and \leref{2.1}
that $A\ot B$ is faithfully flat as a right
$A\square_H B$-module.
\end{proof}

\section{Application to Hopf bimodules}\selabel{4}
Now let $A$ and $B$ be right $H$-comodule algebras.
A two-sided
relative Hopf module is a $k$-module with a left $A$-action, a right $B$-action,
and a right $H$-coaction, such that
$$\rho(amb)= a_{[0]}m_{[0]}b_{[0]}\ot a_{[1]}m_{[1]}b_{[1]},$$
for all $a\in A$, $b\in B$ and $m\in M$.
${}_A\Mm_B^H$ is the category of two-sided relative Hopf modules with
$k$-module maps that are $A$-linear, $B$-linear and $H$-colinear.\\

 $B^{{\rm op}}$ is a
left $H$-comodule algebra, with left coaction $\lambda$ given by
$\lambda(b)=S^{-1}(b_{[1]})\ot b_{[0]}$.
We can then apply the above results to $A$ and $B^{\rm op}$. In particular,
$A\ot B^{\rm op}$ is a right $H\ot H^{\rm cop}$-comodule algebra.

\begin{lemma}\lelabel{2.4}
Let $A$ and $B$ be right $H$-comodule algebras. Then the
Doi-Hopf modules category
${}_{A\ot B^{\rm op}}\Mm(H\ot H^{\rm cop})^H$ is isomorphic to the
category of two-sided relative Hopf modules ${}_A\Mm_B^H$.
\end{lemma}

\begin{proof}
It is well-known that ${}_{A\ot B^{\rm op}}\Mm$ is isomorphic to
the category of bimodules ${}_A\Mm_B$. The isomorphism respects the compatibility of
the action and coaction.
\end{proof}

$A\ot B^{\rm op}$ is a two-sided Hopf module, with coaction
$\rho(a\ot b)=a_{[0]}\ot b_{[0]}\ot a_{[1]}b_{[1]}$.
Furthermore
$(A\ot B^{\rm op})^{{\rm co}H}=A\square_H B^{{\rm op}}$.
Applying \thref{2.3}, we obtain the following Structure Theorem for
two-sided Hopf modules.

\begin{theorem}\thlabel{2.4a}
Let $H$ be a Hopf algebra over the commutative ring $k$, with bijective antipode, and consider two right
$H$-comodule algebras $A$ and $B$.
We have a pair of adjoint functors $(F=A\ot B^{\rm op}\ot_{A\square_H B^{\rm op}}-,\
G=(-)^{{\rm co}H})$ between the categories ${}_{A\square_H B^{\rm op}}\Mm$
and ${}_A\Mm_B^H$. If $A$ is a right faithfully
flat $H$-Galois extension, then $(F,G)$ is a pair of inverse equivalences.
\end{theorem}

\begin{remark} \label{r:equiv} Assume that $A$  (resp. $B$) is a right
(resp. left) faithfully
flat $H$-Galois extension. The proof of \thref{2.3} shows that
via appropriate transport of structure, the functors
$$(A\otimes B^{\rm op})\otimes_{A\square B^{\rm op}}-,\ \
A\otimes_{A^{{\rm co}H}}-,\ \ -\otimes_{B^{{\rm
co}H}}B:{}_{A\square B^{\rm op}}\mathcal{M}
\to{}_A\mathcal{M}^H_B$$ are naturally isomorphic equivalences of
categories. It follows immediately that we may define the functors
$$-\otimes_{A^{{\rm co}H}}-: \ _{B\square A^{\rm op}}\mathcal{M}\times \
_A\mathcal{M}_C\to{}_B\mathcal{M}_C,$$
$$-\otimes_{A^{{\rm co}H}}-: \ _B\mathcal{M}_A\times{}_{A\square
C^{\rm op}}\mathcal{M}\to \ _B\mathcal{M}_C.$$
\end{remark}

\begin{proposition}\prlabel{3.1}
Let $A,B,C$ be right $H$-comodule algebras. If $M\in {}_A\Mm_B^H$
and $N\in {}_B\Mm_C^H$, then $M\ot_B N\in  {}_A\Mm_C^H$. If $A$ and $B$ are right
faithfully flat $H$-Galois extensions, then the map
$$f:\ M^{{\rm co}H}\ot_{B^{{\rm co}H}} N^{{\rm co}H}\to (M\ot_B N)^{{\rm co}H},~~
f(m\ot n)=m\ot n,$$
is an isomorphism. Consequently $M^{{\rm co}H}\ot_{B^{{\rm co}H}} N^{{\rm co}H}$
is a left $A\square_H C^{\rm op}$-module.
\end{proposition}

\begin{proof}
It is clear that $M\ot_B N$ is an $(A,C)$-bimodule. A right $H$-coaction on
$M\ot_B N$ is defined as follows:
$$\rho(m\ot_B n)=m_{[0]}\ot_B n_{[0]}\ot m_{[1]}n_{[1]}.$$
It is easy to show that $\rho$ is well-defined, and that this coaction makes
$M\ot_B N\in {}_A\Mm_C^H$.\\
By restriction of scalars, $M\in {}_A\Mm^H$ and $N\in {}_B\Mm^H$.
It follows from \thref{1.1} that
$$\varepsilon_{1,M}:\ A\ot_{A^{{\rm co}H}} M^{{\rm co}H}\to M~~{\rm and}~~
\varepsilon_{1,N}:\ B\ot_{B^{{\rm co}H}} N^{{\rm co}H}\to N$$
are isomorphisms. Let $g$ be the composition of the maps
$$\varepsilon_{1,M}\ot_{B^{{\rm co}H}} N^{{\rm co}H}:\
A\ot_{A^{{\rm co}H}}M^{{\rm co}H}\ot_{B^{{\rm co}H}}N^{{\rm co}H}
\to M\ot_{B^{{\rm co}H}}N^{{\rm co}H} $$
and
$$M\ot_B\varepsilon_{1,N}:\ M\ot_{B^{{\rm co}H}}N^{{\rm co}H}\cong M\ot_BB\ot_{B^{{\rm co}H}}N^{{\rm co}H}
\to M\ot_BN.$$
$g$ is bijective, and is given by the formula
$$g(a\ot m\ot n)=am\ot_B n,$$
for $a\in A$, $m\in M^{{\rm co}H}$ and $n\in N^{{\rm co}H}$. It is clear that $g$ is
left $H$-linear. $g$ is also right $H$-colinear, since
$g(a_{[0]}\ot m\ot n)\ot a_{[1]}=a_{[0]} m\ot n\ot a_{[1]}=\rho(am\ot_B n)$, and it follows
that $g$ is an isomorphism in $ {}_A\Mm^H$, and, by \thref{1.1} that 
$$g^{{\rm co}H}:\ \bigl(A\ot_{A^{{\rm co}H}}M^{{\rm co}H}\ot_{B^{{\rm co}H}}N^{{\rm co}H}\bigr)^{{\rm co}H}
\to (M\ot_BN)^{{\rm co}H}$$
is an isomorphism. The map $f$ is an isomorphism since it is the composition of $g^{{\rm co}H}$ and the
isomorphism
$$\eta_{1,M^{{\rm co}H}\ot_{B^{{\rm co}H}}N^{{\rm co}H}}:\
M^{{\rm co}H}\ot_{B^{{\rm co}H}}N^{{\rm co}H}\to \bigl(A\ot_{A^{{\rm co}H}}M^{{\rm co}H}\ot_{B^{{\rm co}H}}N^{{\rm co}H}\bigr)^{{\rm co}H}.$$
Finally, the left $A\sq_HC^{\rm op}$-action on $(M\ot_BN)^{{\rm co}H}$ can be transported
using $f$ to $M^{{\rm co}H}\ot_{B^{{\rm co}H}}N^{{\rm co}H}$.
\end{proof}

In the sequel, we will use the adjoint pair of functors $(F,G)$, with unit $\eta$ and counit
$\varepsilon$ introduced after \leref{2.1}, in the cases where the algebras involved are
respectively $A$ and $B^{\rm op}$, $A$ and $C^{\rm op}$ and $B$ and $C^{\rm op}$.
If $A$ and $B$ are right faithfully flat $H$-Galois extensions, then these three adjunctions
are pairs of inverse equivalences, by \thref{2.4a}. We will use the same notation $(F,G)$
and $(\eta,\varepsilon)$ for the three adjunctions, no confusion will arise from this.\\
Take $M_1\in {}_{A\square_H B^{\rm op}}\Mm$ and $N_1\in {}_{B\square_H C^{\rm op}}\Mm$,
and denote
$$M= (A\ot B^{{\rm op}})\ot_{A\square_H B^{\rm op}}M_1\in {}_A\Mm_B^H;$$
$$N= (B\ot C^{{\rm op}})\ot_{B\square_H C^{\rm op}}N_1\in {}_B\Mm_C^H.$$
Using \thref{2.4a} and \prref{3.1}, we find isomorphisms
$$M_1\ot_{B^{{\rm co}H}}N_1\cong M^{{\rm co}H}\ot_{B^{{\rm co}H}} N^{{\rm co}H}
\cong (M\ot_B N)^{{\rm co}H}\in {}_{A\sq_H C^{\rm op}}\Mm.$$
Transporting structure, we find that $M_1\ot_{B^{{\rm co}H}}N_1\in {}_{A\sq_H C^{\rm op}}\Mm$,
and we have a functor
$$-\ot_{B^{{\rm co}H}}-:\ {}_{A\sq_H B^{\rm op}}\Mm\times {}_{B\sq_H C^{\rm op}}\Mm
\to {}_{A\sq_H C^{\rm op}}\Mm.$$

\begin{corollary}\colabel{3.2}
Let $A,B,C$ be right $H$-comodule algebras, and
assume that $A$ and $B$ are right faithfully flat $H$-Galois extensions. Take
$M_1\in {}_{A\square_H B^{\rm op}}\Mm$ and $N_1\in {}_{B\square_H C^{\rm op}}\Mm$.
With notation as above, we have that
$M_1\ot_{B^{{\rm co}H}}N_1\in {}_{A\square_H C^{\rm op}}\Mm$,
and we have an isomorphism 
$$h:\  (A\ot C^{{\rm op}})\ot_{A\square_H C^{\rm op}}
(M_1\ot_{B^{{\rm co}H}}N_1)\to M\ot_B N$$
in ${}_A\Mm_C^H$. This isomorphism is natural in $M_1$ and $N_1$.
\end{corollary}

For later use, we observe that the naturality of $h$ means the following. Let
$\mu_1:\ M_1\to M'_1$ and $\nu_1:\ N_1\to N'_1$ be morphisms in respectively
${}_{A\sq_H B^{\rm op}}\Mm$ and ${}_{B\sq_H C^{\rm op}}\Mm$, and let
$\mu=F(\mu_1)$, $\nu=F(\nu_1)$. Then $\mu_1\ot_{B^{{\rm co}H}}\nu_1$ is a morphism
in ${}_{A\sq_H C^{\rm op}}\Mm$, and the following diagram commutes
\begin{equation}\eqlabel{3.2.2}
\xymatrix{
F(M_1\ot_{B^{{\rm co}H}}N_1)
\ar[rr]^{F (\mu_1\ot \nu_1)}\ar[d]^{h}&&
F(M'_1\ot_{B^{{\rm co}H}}N'_1)\ar[d]^{h}\\
M\ot_BN\ar[rr]^{\mu\ot \nu}&& M'\ot_BN'}
\end{equation}

From now on, let $H$ be a projective Hopf algebra (this condition is always fulfilled
if $k$ is a field); let $A$ be a right $H$-comodule algebra,
and $M,N\in {}_A\Mm^H$. Then the map
$$\nu:\ {}_A\Hom(M,N)\ot H\to {}_A\Hom(M,N\ot H),~~\nu(f\ot h)(m)=f(m)\ot h$$
is injective (see for example \cite[Prop. II.4.2, p. AII.75]{Bourbaki}).
A direct computation shows that the map 
$\tilde{\rho}:\ {}_A\Hom(M,N)\to {}_A\Hom(M,N\ot H)$ defined by
$$\tilde{\rho}(f)(m)=f(m_{[0]})_{[0]}\ot S^{-1}(m_{[1]})f(m_{[0]})_{[1]}$$
is left $A$-linear. Let ${}_A\HOM(M,N)$ be the $k$-submodule of
${}_A\Hom(M,N)$ consisting of the maps $f$ for which
$\tilde{\rho}(f)$
factorizes through ${}_A\Hom(M,N)$, or, equivalently, for which there
exists $f_{[0]}\ot f_{[1]}\in {}_A\Hom(M,N)\ot H$ such that
\begin{equation}\eqlabel{3.2.1}
f_{[0]}(m)\ot f_{[1]}=f(m_{[0]})_{[0]}\ot S^{-1}(m_{[1]})f(m_{[0]})_{[1]},
\end{equation}
for all $m\in M$. It follows from the injectivity of $\nu$ that $f_{[0]}\ot f_{[1]}$
is unique if it exists. ${}_A\HOM(M,N)$ is called the rational part of
${}_A\Hom(M,N)$. If $H$ is finitely generated and projective, then
$\nu$ is bijective, and ${}_A\HOM(M,N)={}_A\Hom(M,N)$. We have a map
$$\rho=\nu^{-1}\circ \tilde{\rho}:\ {}_A\HOM(M,N)\to {}_A\Hom(M,N)\ot H,~~
\rho(f)=f_{[0]}\ot f_{[1]}.$$

\begin{proposition}\prlabel{3.2x}
Let $H$ be a projective Hopf algebra, $A$ a right $H$-comodule algebra,
and $M,N\in {}_A\Mm^H$. Then $({}_A\HOM(M,N),\rho)$ is a right
$H$-comodule.
\end{proposition}

\begin{proof}
$N\ot H\in {}_A\Mm^H$ under the diagonal coaction.
We know that five of the six faces of the following diagram, namely all faces except the
top one, commute.
$$\xymatrix{
\sf{}_A\HOM(M,N)\ar[rr]^{\rho}\ar[dd]^{\subset}\ar[dr]^{\rho}&&\sf{}_A\Hom(M,N)\ot H\ar'[d][dd]^{\nu}
\ar[dr]^{\tilde{\rho}\ot H}&\\
&\sf{}_A\Hom(M,N)\ot H\ar[dd]^(.7){\nu}\ar[rr]^(.35){(\nu\circ H)\circ (id\ot \Delta)}
&&\sf{}_A\Hom(M,N\ot H)\ot H\ar[dd]^{\nu}\\
\sf{}_A\Hom(M,N)\ar'[r][rr]^(.3){\tilde{\rho}}\ar[dr]^{\tilde{\rho}}&&\sf{}_A\Hom(M,N\ot H)\ar[dr]^{\tilde{\rho}}&\\
&\sf{}_A\Hom(M,N)\ot H)\ar[rr]^{{}_A\Hom(M,N\ot \Delta)}
&&\sf{}_A\Hom(M,N\ot H\ot H)}$$
This implies that the top face also commutes; this means that, for all
$f\in {}_A\HOM(M,N)$,
$$(\nu\ot H)(f_{[0]}\ot \Delta(f_{[1]}))=\tilde{\rho}(f_{[0]})\ot f_{[1]},$$
and therefore $f_{[0]}\ot f_{[1]}\in {}_A\HOM(M,N)\ot H$. We then also have that
$$(\nu\ot H)(f_{[0]}\ot \Delta(f_{[1]}))=(\nu\ot H)(\rho(f_{[0]})\ot f_{[1]}),$$
and, since $\nu\ot H$ is injective,
$$f_{[0]}\ot \Delta(f_{[1]})=\rho(f_{[0]})\ot f_{[1]}.$$
We therefore have shown that $\rho:\ {}_A\HOM(M,N)\to {}_A\HOM(M,N)\ot H$
is a coassociative map. Finally, it follows immediately from \equref{3.2.1}
that $\varepsilon(f_{[1]})f_{[0]}=f$, for all $f\in  {}_A\HOM(M,N)$.
\end{proof}

An alternative description of ${}_A\HOM(M,N)$ is the following:
${}_A\Hom(M,N)$ is a left $H^*$-module,
with action (see \cite[6.5]{DNR} in the case where $k$ is a field):
$$(h^*\cdot f)(m)=\lan h^*, S^{-1}(m_{[1]})f(m_{[0]})_{[1]}\ran f(m_{[0]})_{[0]}.$$
${}_A\HOM(M,N)$ is then the subspace of ${}_A\Hom(M,N)$
consisting of left $A$-linear $f:\ M\to N$ for which there exists
a (unique) $f_{[0]}\ot f_{[1]}\in {}_A\Hom(M,N)\ot H$ such that
$$(h^*\cdot f)(m)=\lan h^*,f_{[1]}\ran f_{[0]}(m).$$

\begin{proposition}\prlabel{3.2c}
Let $A$ be a right $H$-comodule algebra, with $H$ a projective Hopf algebra,
and $M,N\in {}_A\Mm^H$. If $M$ is finitely generated projective as a left $A$-module, then
${}_A\HOM(M,N)$ coincides with ${}_A\Hom(M,N)$. For $f\in {}_A\HOM(M,N)$,
we have
\begin{equation}\eqlabel{3.2c.1}
\rho(f)=\sum_i m_i^*\cdot f(m_{i[0]})_{[0]}\ot S^{-1}(m_{i[1]})f(m_{i[0]})_{[1]},
\end{equation}
where $\sum_i m_i^*\ot_A m$ is a finite dual basis of $M\in {}_A\Mm$.
\end{proposition}

\begin{proof}
We used the following notation: for $m^*\in {}_A\Hom(M,A)$, and $n\in N$, $m^*\cdot n
\in {}_A\Hom(M,N)$ is defined by
$$(m^*\cdot n)(m)=m^*(m)n.$$
For every $m\in M$, we have that $m=\sum_i m_i^*(m)m_i$, hence
\begin{equation}\eqlabel{3.2c.2}
\rho(m)=\sum_i m_i^*(m)_{[0]}m_{i[0]}\ot m_i^*(m)_{[1]}m_{i[1]}.
\end{equation}
We then compute that
\begin{eqnarray*}
&&\hspace*{-1cm}
f(m_{[0]})_{[0]}\ot S^{-1}(m_{[1]})f(m_{[0]})_{[1]}\\
&\equal{\equref{3.2c.2}}&
\sum_i f(m_i^*(m)_{[0]}m_{i[0]})_{[0]} \ot S^{-1}(m_i^*(m)_{[1]}m_{i[1]})f(m_i^*(m)_{[0]}m_{i[0]})_{[1]}\\
&=&
\sum_i m_i^*(m)_{[0]}f(m_{i[0]})_{[0]} \ot S^{-1}(m_{i[1]})S^{-1} (m_i^*(m)_{[2]})
m_i^*(m)_{[1]}f(m_{i[0]})_{[1]}\\
&=& \sum_i m_i^*(m)f(m_{i[0]})_{[0]}\ot S^{-1}(m_{i[1]})
(m_{i[0]})_{[1]}\\
&=& \sum_i m_i^*\cdot f(m_{i[0]})_{[0]}(m)\ot S^{-1}(m_{i[1]})f(m_{i[0]})_{[1]},
\end{eqnarray*}
and \equref{3.2c.1} follows from \equref{3.2.1}.
\end{proof}

\begin{proposition}\prlabel{3.3}
Let $H$ be a projective Hopf algebra, and $A,B,C$ right $H$-comodule
algebras. If $M\in {}_A\Mm_B^H$ and $N\in {}_A\Mm_C^H$, then
$${}_A\HOM(M,N)\in {}_B\Mm_C^H.$$
We have a map
$$\beta:\ {}_A\HOM(M,N)^{{\rm co}H}\to {}_{A^{{\rm co}H}}\Hom(M^{{\rm co}H},N^{{\rm co}H}).$$
If $A$ is a right faithfully flat $H$-Galois extension, then $\beta$ is an isomorphism of left $B\square C^{\rm op}$-modules.
\end{proposition}

\begin{proof}
We consider the following $(B,C)$-bimodule structure on ${}_A\Hom(M,N)$:
$$(b\cdot f\cdot c)(m)=f(mb)c.$$
It is clear that $b\cdot f\cdot c$ is then left $A$-linear.
Take $f\in {}_A\HOM(M,N)$; in order to show that $b\cdot f\cdot c\in {}_A\HOM(M,N)$,
it suffices to show that $b_{[0]}\cdot f_{[0]}\cdot c_{[0]}\ot b_{[1]}f_{[1]}c_{[1]}$
satisfies \equref{3.2.1}. This can be seen as follows: for all $m\in M$, we have
\begin{eqnarray*}
&&\hspace*{-2cm}
(b_{[0]}\cdot f_{[0]}\cdot c_{[0]})(m)\ot b_{[1]}f_{[1]}c_{[1]}\\
&=& f_{[0]}(mb_{[0]})c_{[0]}\ot b_{[1]}f_{[1]}c_{[1]}\\
&=& f\bigl((mb_{[0]})_{[0]}\bigr)_{[0]}c_{[0]}\ot b_{[1]}S^{-1}\bigl((mb_{[0]})_{[1]}\bigr)
 f\bigl((mb_{[0]})_{[0]}\bigr)_{[1]}c_{[1]}\\
&=& f(m_{[0]}b_{[0]})_{[0]}c_{[0]}\ot b_{[2]}S^{-1}(b_{[1]})S^{-1}(m_{[1]})
f(m_{[0]}b_{[0]})_{[1]}c_{[1]}\\
&=& f(m_{[0]}b)_{[0]}c_{[0]}\ot S^{-1}(m_{[1]})f(m_{[0]}b)_{[1]} c_{[1]}\\
&=& \bigl(f(m_{[0]}b)c\bigr)_{[0]}\ot S^{-1}(m_{[1]})  \bigl(f(m_{[0]}b)c\bigr)_{[1]}\\
&=& (b\cdot f\cdot c)(m_{[0]})_{[0]}\ot S^{-1}(m_{[1]}) (b\cdot f\cdot c)(m_{[0]})_{[1]},
\end{eqnarray*}
as needed. This shows also that $\rho(b\cdot f\cdot c)=
b_{[0]}\cdot f_{[0]}\cdot c_{[0]}\ot b_{[1]}f_{[1]}c_{[1]}$, hence that
${}_A\HOM(M,N)\in {}_B\Mm_C^H$.\\
Now take $ f\in {}_A\HOM(M,N)^{{\rm co}H}$. Then $\rho(f)=f_{[0]}\ot f_{[1]}$, so
$$f(m)\ot 1=f(m_{[0]})_{[0]}\ot S^{-1}(m_{[1]})f(m_{[0]})_{[1]},$$
for all $m\in M$. If $m\in M^{{\rm co}H}$, then it follows that $f(m)\ot 1=\rho(f(m))$,
so $f(m)\in N^{{\rm co}H}$. Thus $f$ restricts to a map $\beta(f)=f^{{\rm co}H}:\
M^{{\rm co}H}\to N^{{\rm co}H}$. Using the fact that $f$ is left $A$-linear, we see
that the diagram
$$\diagram
A\ot_{A^{{\rm co}H}} M^{{\rm co}H}\rrto^{A\ot f^{{\rm
co}H}}\dto_{\varepsilon_{1,M}} &&
A\ot_{A^{{\rm co}H}} N^{{\rm co}H}\dto^{\varepsilon_{1,N}}\\
M\rrto^{f}&& N
\enddiagram$$
commutes. If $A$ is right faithfully flat $H$-Galois, then we define the
inverse of $\beta$ as follows:
$$\beta^{-1}(g)=\varepsilon_{1,N}\circ (A\ot g)\circ \varepsilon_{1,M}^{-1}.$$
\end{proof}

Combining \prref{3.3} with \thref{2.4a}, we obtain the following
result.

\begin{corollary}\colabel{3.4}
Let $A,B,C$ be right $H$-comodule algebras, and assume that $A$ and $B$ are right faithfully flat
$H$-Galois extensions.
Let $M_1\in {}_{A\square B^{\rm op}}\Mm$ and $N_1\in {}_{A\square C^{\rm op}}\Mm$,
and consider
$$M= (A\ot B^{\rm op})\ot_{A\square B^{\rm op}}M_1\in {}_A\Mm_B^H,$$
$$N= (A\ot C^{\rm op})\ot_{A\square C^{\rm op}}N_1\in {}_A\Mm_C^H.$$
Then
$${}_{A^{{\rm co}H}}\Hom(M_1,N_1)
\cong {}_A\HOM(M,N)^{{\rm co}H}\in {}_{B\square C^{\rm op}}
\Mm$$
and
$${}_A\HOM(M,N)\cong (B\ot C^{\rm op})\ot_{B\square C^{\rm op}}
{}_{A^{{\rm co}H}}\Hom(M_1,N_1).$$
\end{corollary}

\begin{proposition}\prlabel{3.5}
Let $A,B,C$ be right $H$-comodule algebras, and consider $M\in {}_A\Mm_B^H$,
$N\in {}_A\Mm_C^H$. Then the evaluation map
$$\varphi:\ M\ot_B {}_A\HOM(M,N)\to N,~~\varphi(m\ot_B f)=f(m)$$
is a morphism in ${}_A\Mm_C^H$.\\
If $A$ and $B$ are right faithfully flat $H$-Galois extensions, then the evaluation map
$$M^{{\rm co}H}\ot_{B^{{\rm co}H}} {}_{A^{{\rm co}H}}\Hom(M^{{\rm co}H},N^{{\rm co}H})
\to N^{{\rm co}H}$$
is left $A\square_H C^{\rm op}$-linear.
\end{proposition}

\begin{proof}
We first show that $\varphi$ is right $H$-colinear.
\begin{eqnarray*}
&&\hspace*{-2cm}
(\varphi\ot H)(\rho(m\ot f))=(\varphi\ot H)(m_{[0]}\ot_B f_{[0]}\ot m_{[1]}f_{[1]}) \\
&=& f_{[0]}(m_{[0]})\ot m_{[1]}f_{[1]}\\
&\equal{\equref{3.2.1}}&f(m_{[0]})_{[0]}\ot m_{[2]}S^{-1}(m_{[1]})f(m_{[0]})_{[1]}\\
&=& \rho(f(m))=\rho(\varphi(m\ot_B f)).
\end{eqnarray*}
$\varphi$ is left $A$-linear and right $C$-linear since
$$\varphi(am\ot f\cdot c)=(f\cdot c)(am)=f(am)c=af(m)c=a\varphi(m\ot_B f)c.$$
The composition
\begin{eqnarray*}
&&\hspace*{-8mm} M^{{\rm co}H}\ot_{B^{{\rm co}H}} {}_{A^{{\rm
co}H}}\Hom(M^{{\rm co}H},N^{{\rm co}H}) 
\stackrel{M^{{\rm
co}H}\ot\beta}{\longrightarrow}
M^{{\rm co}H}\ot_{B^{{\rm co}H}} {}_A\HOM(M,N)^{{\rm co}H}\\
&&\hspace*{1cm}
\stackrel{f}{\longrightarrow} (M\ot_B {}_A\HOM(M,N))^{{\rm co}H}
\stackrel{\varphi^{{\rm co}H}}{\longrightarrow} N^{{\rm co}H}
\end{eqnarray*}
is the required evaluation map. If $A$ is right faithfully flat $H$-Galois, then
$\beta$ is an isomorphism of $B\square_H C^{\rm op}$-modules, by \prref{3.3},
and then $M^{{\rm co}H}\ot\beta$ is an isomorphism of $A\square_H C^{\rm op}$-modules,
by \coref{3.2}.
If $B$ is right faithfully flat $H$-Galois, then
$f$ is an isomorphism of $A\square_H C^{\rm op}$-modules, by \prref{3.1}.
$\varphi$ is a morphism in ${}_A\Mm_C^H$, hence
$\varphi^{{\rm co}H}$ is left $A\square_H C^{\rm op}$-linear, since $(-)^{{\rm co}H}$
is a functor from ${}_A\Mm_C^H$ to ${}_{A\square_H C^{\rm op}}\Mm$.
\end{proof}

\begin{proposition}\prlabel{3.6}
Let $A$ be a right $H$-comodule algebra, and $M\in {}_A\Mm^H$. Then ${}_A\END(M)^{\rm op}$
is a right $H$-comodule algebra.
\end{proposition}

\begin{proof}
Applying \prref{3.3} (with $M=N$, $B=C=k$), we see that ${}_A\END(M)$ is a right $H$-comodule.
We have to show the compatibility relation
\begin{equation}\eqlabel{3.6.1}
\rho(g\circ f)=g_{[0]}\circ f_{[0]}\ot f_{[1]}g_{[1]},
\end{equation}
for all $f,g\in {}_A\END(M)$. To this end, it suffices to show that the right hand side of \equref{3.6.1}
satisfies \equref{3.2.1}. Indeed, for all $m\in M$, we have
\begin{eqnarray*}
&&\hspace*{-13mm}
(g_{[0]}\circ f_{[0]})(m)\ot f_{[1]}g_{[1]}\\
&\equal{\equref{3.2.1}}&
g_{[0]}\bigl(f(m_{[0]})_{[0]}\bigr)\ot S^{-1}(m_{[1]})f(m_{[0]})_{[1]}g_{[1]}\\
&\equal{\equref{3.2.1}}&
g\bigl(f(m_{[0]})_{[0]}\bigr)_{[0]}\ot S^{-1}(m_{[1]})f(m_{[0]})_{[2]}
S^{-1}\bigl(f(m_{[0]})_{[1]}\bigr)g\bigl(f(m_{[0]})_{[0]}\bigr)_{[1]}\\
&=&g\bigl(f(m_{[0]})\bigr)_{[0]}\ot S^{-1}(m_{[1]})g\bigl(f(m_{[0]})\bigr)_{[1]}\\
&=&(g\circ f)(m_{[0]})_{[0]}\ot S^{-1}(m_{[1]}) (g\circ f)(m_{[0]})_{[1]}.
\end{eqnarray*}
\end{proof}

\begin{proposition}\prlabel{3.7}
Let $A,B$ be right $H$-comodule algebras, and consider
$M\in {}_A\Mm_B^H$. Then the map
$$\psi:\ B\to {}_A\END(M),~~\psi(b)(m)=mb$$
is a morphism in ${}_B\Mm_B^H$.\\
If $A$ is a right faithfully flat $H$-Galois extension, then the map
$$\psi^{{\rm co}H}:\ B^{{\rm co}H}\to {}_A\END(M)^{{\rm co}H}\cong {}_{A^{{\rm co}H}}
\End(M^{{\rm co}H})$$
is left $B\square_H B^{\rm op}$-linear.
\end{proposition}

\begin{proof}
We first show that $\psi$ is right $H$-colinear and well-defined. Indeed,
$$\psi(b)_{[0]}\ot \psi(b)_{[1]}=\psi(b_{[0]})\ot b_{[1]},$$
since
\begin{eqnarray*}
&&\hspace*{-2cm}
\psi(b)(m_{[0]})_{[0]}\ot S^{-1}(m_{[1]}) \psi(b)(m_{[0]})_{[1]}\\
&=& (m_{[0]}b)_{[0]} \ot S^{-1}(m_{[1]}) (m_{[0]}b)_{[1]}\\
&=& m_{[0]}b_{[0]}\ot S^{-1}(m_{[2]})m_{[1]}b_{[1]}\\
&=& mb_{[0]}\ot b_{[1]}=\psi(b_{[0]})(m)\ot b_{[1]}.
\end{eqnarray*}
$\psi$ is left and right $B$-linear since
$$\psi(b'bb'')(m)=mb'bb''=((b'\cdot \psi\cdot b'')(b))(m),$$ for all $b,b',b''\in B$
and $m\in M$. The second statement then follows immediately from \coref{3.4}.
\end{proof}

\begin{remark}\relabel{3.7b}
The map $\psi$ in \prref{3.7} is also a morphism of right $H$-comodule algebras
between $B$ and ${}_A\END(M)^{\rm op}$.
\end{remark}

\begin{proposition}\prlabel{3.8}
Let $A,B,C$ be right $H$-comodule algebras, and consider
$M\in {}_A\Mm_B^H$, $N\in {}_A\Mm_C^H$. Then the map
$$\mu:\ {}_A\HOM(M,A)\ot_A N\to {}_A\HOM(M,N),~~
\mu(f\ot n)(m)=f(m)n$$
is a morphism in ${}_B\Mm_C^H$. If $A$ is a right faithfully flat $H$-Galois extension, then
the map
\begin{eqnarray*}
&&\hspace*{-1cm}
\mu^{{\rm co}H}:\ {}_{A^{{\rm co}H}}\Hom(M^{{\rm co}H},A^{{\rm co}H})\ot_{A^{{\rm co}H}}
{N^{{\rm co}H}}
\cong ({}_A\HOM(M,A)\ot_A N)^{{\rm co}H}\\
&&\hspace*{1cm}\to
{}_{A^{{\rm co}H}}\Hom(M^{{\rm co}H},N^{{\rm co}H})\cong {}_A\HOM(M,N)^{{\rm co}H}
\end{eqnarray*}
is left $B\square_H C^{\rm op}$-linear.
\end{proposition}

\begin{proof}
In order to prove that $\mu$ is right $H$-colinear, we have to show that
$$\rho(\mu(f\ot n))=\mu(f_{[0]}\ot n_{[0]})\ot f_{[1]}n_{[1]}.$$
It suffices to compute that
\begin{eqnarray*}
&&\hspace*{-2cm}
\mu(f_{[0]}\ot n_{[0]})(m)\ot f_{[1]}n_{[1]}
\equal{\equref{3.2.1}} f(m_{[0]})_{[0]}n_{[0]}\ot S^{-1}(m_{[1]})f(m_{[0]})_{[1]}n_{[1]}\\
&=& (f(m_{[0]})n)_{[0]}\ot S^{-1}(m_{[1]})(f(m_{[0]})n)_{[1]}\\
&=& (\mu(f\ot n)(m_{[0]}))_{[0]}\ot S^{-1}(m_{[1]})(\mu(f\ot n)(m_{[0]}))_{[1]}.
\end{eqnarray*}

Finally, $\mu$ is left $B$-linear and right $C$-linear, since
$$(\mu(b\cdot f\ot nc))(m)=f(mb)nc=\mu(f\ot n)(mb)c
=(b\cdot \mu(f\ot n)\cdot c)(m).$$
\end{proof}

\section{Morita equivalences}\selabel{5}
In this section, we study Morita equivalences induced by two-sided relative
Hopf modules.

\begin{definition}\delabel{5.1}
Let $A$ and $B$ be right $H$-comodule algebras. An $H$-Morita context
connecting $A$ and $B$ is a Morita context $(A,B,M,N,\alpha,\beta)$ such that
$M\in {}_A\Mm_B^H$, $N\in {}_B\Mm_A^H$,  $\alpha:\ M\ot_B N\to A$
is a morphism in ${}_A\Mm_A^H$ and $\beta:\ N\ot_A M\to B$ is a morphism in
${}_B\Mm_B^H$.
\end{definition}

A morphism between two $H$-Morita contexts $(A,B,M,N,\alpha,\beta)$
and $(A',B',$ $M',N',\alpha',\beta')$ is defined in the obvious way: it consists of a 
fourtuple $(\kappa,\lambda,$ $\mu,\nu)$, where $\kappa:\ A\to A'$ and $\lambda:\
B\to B'$ are $H$-comodule algebra maps, $\mu:\ M\to M'$ is a morphism in
${}_A\Mm_B^H$ and $\nu:\ N\to N'$ is a morphism in ${}_B\Mm_A^H$
such that $\kappa\circ\alpha=\alpha'\circ (\mu\ot\nu)$ and $\lambda\circ\beta
=\beta'\circ (\nu\ot\mu)$. $\Morita^H(A,B)$ will be the subcategory of the category of
$H$-Morita contexts, consisting of $H$-Morita contexts connecting $A$ and $B$,
and morphisms with the identity of $A$ and $B$ as the underlying algebra maps.

\begin{proposition}\prlabel{5.2}
Let $(A,B,M,N,\alpha,\beta)$ be a strict $H$-Morita context. Then
we have a pair of inverse equivalences $(M\ot_B-, N\ot_A-)$
between the categories ${}_A\Mm^H$ and ${}_B\Mm^H$.
\end{proposition}

\begin{proof}
Let $P\in {}_B\Mm^H$. Then $M\ot_B P\in {}_A\Mm^H$, with right $H$-action
$$\rho(m\ot_B p)=m_{[0]}\ot_B p_{[0]}\ot m_{[1]}p_{[1]}.$$ The rest of the proof is
straightforward.
\end{proof}

We will now give an $H$-comodule version of \cite[Prop. 4.2.1]{Bass}.

\begin{example}\exlabel{5.3}
Let $A$ be a right $H$-comodule algebra, and $M\in {}_A\Mm^H$. Then
$B={}_A\END(M)^{\rm op}$ is also a right $H$-module algebra, by \prref{3.6}.
Then $M\in {}_A\Mm^H_B$, with right $B$-action given by $m\cdot f=f(m)$, for all
$f\in B$ and $m\in M$. Indeed, $(m\cdot f)\cdot g=m\cdot (g\circ f)$, and
\begin{eqnarray*}
&&\hspace*{-2cm}
m_{[0]}\cdot f_{[0]}\ot m_{[1]}f_{[1]}= f_{[0]}(m_{[0]})\ot m_{[1]}f_{[1]}\\
&\equal{\equref{3.2.1}}& f(m_{[0]})_{[0]}\ot m_{[2]}S^{-1}(m_{[1]})f(m_{[0]})_{[1]}\\
&=&\rho(f(m))=\rho(m\cdot f).
\end{eqnarray*}
It follows from \prref{3.3} that $N={}_A\HOM(M,A)\in {}_B\Mm_A^H$, and from
\prref{3.5} that the map
$$\alpha:\ M\ot_BN\to A,~~\alpha(m\ot n)=n(a)$$
is a morphism in ${}_A\Mm_A^H$. It follows from \prref{3.8} that the map
$$\beta:\ N\ot_A M\to {}_A\END(M),~~\beta(n\ot m)(x)=n(x)m$$
is a morphism in  ${}_B\Mm_B^H$. Straightforward computations then show that
$(A,B,M,N,\alpha,\beta)$ is an $H$-Morita context. We call it the $H$-Morita context
associated to $M\in {}_A\Mm^H$.
\end{example}

\begin{proposition}\prlabel{5.4}
The $H$-Morita context associated to $M\in {}_A\Mm^H$ is strict if and only if $M$ is
a progenerator as a left $A$-module.
\end{proposition}

\begin{proof}
If the Morita context is strict, then $M$ is a left $A$-progenerator by \cite[Theorem III.3.5]{Bass}.
Conversely, if $M$ is a left $A$-progenerator, then $M\in {}_A\Mm$ is finitely generated and
projective, hence ${}_A\Hom(M,X)={}_A\HOM(M,X)$, for all $X\in {}_A\Mm^H$. If we forget
the $H$-comodule structure in the $H$-Morita context, then we obtain the Morita context
associated to $M\in {}_A\Mm$, as in \cite[Prop. II.4.1]{Bass}. By \cite[Prop. II.4.4]{Bass}, this
Morita context is strict.
\end{proof}

\begin{proposition}\prlabel{5.5}
Let $(A,B,M,N,\alpha,\beta)$ be a strict $H$-Morita context. Then $M$ is a left $A$-progenerator,
and the $H$-Morita context is isomorphic to the $H$-Morita context associated to
$M\in {}_A\Mm^H$.
\end{proposition}

\begin{proof}
$M$ is a left $A$-progenerator by \cite[Theorem III.3.5]{Bass}. Then ${}_A\End(M)={}_A\END(M)$,
and by \cite[Theorem II.3.4]{Bass}, $\psi:\ B\to {}_A\END(M)^{\rm op}$, $\psi(b)(m)=mb$
is an isomorphism of $k$-algebras. It is an isomorphism of $H$-comodule algebras, by
\reref{3.7b}. It follows from \cite[Theorem 3.4]{Bass} that
$$\varphi:\ N\to {}_A\HOM(M,A)={}_A\Hom(M,A),~~\varphi(n)(m)=\alpha(m\ot n)$$
is an isomorphism of $(B,A)$-bimodules. We verify that $\varphi$ is $H$-colinear.
For every $n\in N$, we have to show that
\begin{equation}\eqlabel{5.5.1}
\varphi(n_{[0]})\ot n_{[1]}=\varphi(n)_{[0]}\ot \varphi(n)_{[1]}.
\end{equation}
Using the right $H$-colinearity of $\alpha$, we find
\begin{eqnarray*}
&&\hspace*{-15mm}
\alpha(m_{[0]}\ot_B n)_{[0]}\ot S^{-1}(m_{[1]})\alpha(m_{[0]}\ot_B n)_{[1]}\\
&=& \alpha(m_{[0]}\ot_B n_{[0]})\ot S^{-1}(m_{[2]})m_{[1]}n_{[1]}
= \alpha(m\ot_B n_{[0]})\ot n_{[1]},
\end{eqnarray*}
and \equref{5.5.1} follows from \equref{3.2.1}. From classical Morita theory (see \cite{Bass}), we
know that $(A,\psi,M,\varphi)$ is an isomorphism of Morita contexts; since $\psi$ and $\varphi$
are $H$-colinear, it follows that is an isomorphism of $H$-Morita contexts.
\end{proof}

\begin{definition}\delabel{5.6}
Assume that $A$ and $B$ are right faithfully flat $H$-Galois extensions of $A^{{\rm co}H}$ and $B^{{\rm co}H}$.
A $\sq_H$-Morita context between $A^{{\rm co}H}$ and $B^{{\rm co}H}$ is a Morita context
$(A^{{\rm co}H},B^{{\rm co}H},M_1,N_1,\alpha_1,\beta_1)$ such that $M_1$ (resp. $N_1$) is a left
$A\sq_HB^{\rm op}$-module (resp. $B\sq_HA^{\rm op}$-module) and
\begin{itemize}
\item $\alpha_1:\ M_1\ot_{B^{{\rm co}H}} N_1\to A^{{\rm co}H}$ is left $A\sq_HA^{\rm op}$-linear,
\item $\beta_1:\ N_1\ot_{A^{{\rm co}H}} M_1\to B^{{\rm co}H}$ is left $B\sq_HB^{\rm op}$-linear.
\end{itemize}
\end{definition}

A morphism between two $\sq_H$-Morita contexts connecting $A^{{\rm co}H}$ and $B^{{\rm co}H}$,
is a morphism between Morita contexts of the form $(A^{{\rm co}H},B^{{\rm co}H},\mu_1,\nu_1)$,
where $\mu_1$ is left $A\sq_HB^{\rm op}$-linear and $\nu_1$ is left $B\sq_HA^{\rm op}$-linear.
The category of $\sq_H$-Morita contexts connecting $A^{{\rm co}H}$ and $B^{{\rm co}H}$
will be denoted by $\Morita^{\sq_H}(A^{{\rm co}H},B^{{\rm co}H})$.

\begin{theorem}\thlabel{5.7}
Let $A$ and $B$ be right faithfully flat $H$-Galois extensions of $A^{{\rm co}H}$ and $B^{{\rm co}H}$.
Then the categories $\Morita^H(A,B)$ and $\Morita^{\sq_H}(A^{{\rm co}H},$ $B^{{\rm co}H})$
are equivalent. The equivalence functors send strict contexts to strict contexts.
\end{theorem}

\begin{proof}
Let $(A,B,M,N,\alpha,\beta)$ be an $H$-Morita context. It follows from
\thref{2.4a} that $M^{{\rm co}H}
\in{}_{A\square_H B^{\rm op}}\Mm$, and $N^{{\rm co}H}
\in {}_{B\square_H A^{\rm op}}\Mm$.
It follows from \prref{3.1} that we have a left $A\square_H A^{\rm op}$-linear map
$$\alpha_1=\alpha^{{\rm co}H}\circ f:\ M^{{\rm co}H}\ot_{B^{{\rm co}H}} N^{{\rm co}H}
\to (M\ot _B N)^{{\rm co}H}\to A^{{\rm co}H},$$
and
a left $B\square_H B^{\rm op}$-linear isomorphism
$$\beta_1=\beta^{{\rm co}H}\circ f:\ N^{{\rm co}H}\ot_{A^{{\rm co}H}} M^{{\rm co}H}
\to (N\ot _A M)^{{\rm co}H}\to B^{{\rm co}H}.$$
From the description of $f$ in \prref{3.1}, it follows that we have a commutative diagram
of isomorphisms
$$\diagram
M^{{\rm co}H}\ot_{B^{{\rm co}H}}N^{{\rm co}H}\ot_{A^{{\rm co}H}}M^{{\rm co}H}
\rto\dto & (M\ot_BN)^{{\rm co}H}\ot_{A^{{\rm co}H}}M^{{\rm co}H}\dto\\
M^{{\rm co}H}\ot_{B^{{\rm co}H}}(N\ot_A M)^{{\rm co}H}\rto &
(M\ot_BN\ot_A M)^{{\rm co}H}
\enddiagram$$
Now $\alpha\ot_A M= M\ot_B \beta$ implies $(\alpha\ot_A M)^{{\rm
co}H}= (M\ot_B \beta)^{{\rm co}H}$, and it follows that
$$\alpha_1\ot_{A^{{\rm co}H}} M^{{\rm co}H}=M^{{\rm
co}H}\ot_{B^{{\rm co}H}}\beta.$$ In a similar way, we have that
$$\beta_1\ot_{B^{{\rm co}H}} N^{{\rm co}H}=N^{{\rm
co}H}\ot_{A^{{\rm co}H}}\alpha$$ and it follows that $(A^{{\rm
co}H}, B^{{\rm co}H}, M^{{\rm co}H},N^{{\rm
co}H},\alpha_1,\beta_1)$ is a Morita context. If
$(A,B,M,N,\alpha,\beta)$ is strict, then $(A^{{\rm co}H}, B^{{\rm
co}H}, M^{{\rm co}H},$ $N^{{\rm co}H},\alpha_1,\beta_1)$ is also
strict.\\
Conversely, let $(A^{{\rm co}H},B^{{\rm co}H},M_1,N_1,\alpha_1,\beta_1)$
be a $\sq_H$-Morita context. Then
$M=F(M_1)=(A\ot B^{\rm op})\ot_{A\sq_H B^{\rm op}}M_1\in {}_A\Mm_B^H$ and
$N=F(N_1)=(B\ot A^{\rm op})\ot_{B\sq_H A^{\rm op}}N_1\in {}_B\Mm_A^H$.
Also observe that $A\cong F(A^{{\rm co}H})= (A\ot A^{\rm op})\ot_{A\sq_H A^{\rm op}}A^{{\rm co}H}$ and
$B\cong F(B^{{\rm co}H})=(B\ot B^{\rm op})\ot_{B\sq_H B^{\rm op}}B^{{\rm co}H}$. We define
$\alpha:\ M\ot_B N\to A$ and $\beta:\ N\ot_A M\to N$ by the commutativity of the following
two diagrams, where the isomorphisms $h$ are defined as in \coref{3.2}.
\begin{equation}\eqlabel{5.7.1}
\xymatrix{
~~F(M_1\ot_{B^{{\rm co}H}}N_1)~~
\ar[r]^-{F (\alpha_1)}\ar[d]^{h}&
~~F(A^{{\rm co}H})~~\ar[d]^{\cong}\\
M\ot_BN\ar[r]^-{\alpha}& A}
\end{equation}
\begin{equation}\eqlabel{5.7.2}
\xymatrix{
~~F(N_1\ot_{B^{{\rm co}H}}M_1)~~
\ar[r]^-{F (\beta_1)}\ar[d]^{h}&
~~F(B^{{\rm co}H})~~\ar[d]^{\cong}\\
N\ot_AM\ar[r]^-{\beta}& B}
\end{equation}
It is clear that $\alpha\in {}_A\Mm_B^H$ and $\beta\in {}_B\Mm_A^H$. We claim that
$(A,B,M,N,\alpha,\beta)$ is an $H$-Morita context. To this end, consider the following
diagram
$$\xymatrix{
M\ot_B N\ot_A M \ar[rr]^-{M\ot\beta}\ar[d]^{h^{-1}}&&
M\ot_BB\ar[d]^{h^{-1}}\\
F(M_1\ot_{B^{{\rm co}H}}N_1\ot_{A^{{\rm co}H}}M_1)\hbox{\hspace*{5mm}}
\ar[rr]^-{F(M_1\ot \beta_1)}\ar[d]^{=}&&
\hbox{\hspace*{5mm}}F(M_1\ot_{B^{{\rm co}H}}B^{{\rm co}H})\ar[d]^{\cong}\\
F(M_1\ot_{B^{{\rm co}H}}N_1\ot_{A^{{\rm co}H}}M_1)\hbox{\hspace*{5mm}}
\ar[rr]^-{F(\alpha_1\ot N_1)}\ar[d]^{h}&&
\hbox{\hspace*{5mm}}F(A^{{\rm co}H}\ot_A^{{\rm co}H}N_1)\ar[d]^{h}\\
M\ot_B N\ot_A M\ar[rr]^-{\alpha\ot_A M}&&A\ot_A M}$$
The top square and the bottom square commute by the definition of $\alpha$ and $\beta$,
and because of the naturality of $h$ (see \equref{3.2.2}). The square in the middle
commutes because $(A^{{\rm co}H},B^{{\rm co}H},M_1,N_1,\alpha_1,\beta_1)$
is a Morita context. So the whole diagram commutes. The composition of the left
vertical morphisms is the identity of $M\ot_BN\ot_A M$, and the composition of the
right vertical morphisms is the natural isomorphism $M\ot_B B\cong A\ot_AM$. So
it follows that the diagram
$$\xymatrix{
M\ot_B N\ot_A M \ar[rr]^-{M\ot\beta}\ar[d]^{\alpha\ot M}&&M\ot_BB\ar[d]^{\cong}\\
A\ot_AM\ar[rr]^-{\cong}&&M}$$
commutes. The commutativity of the second diagram in the definition of a Morita context
is proved in a similar way.
\end{proof}

Recall that $M\in {}_A\Mm$ is a progenerator if and only if $A$ and $M$ are mutually
direct summands of finite direct sums of copies of the other. Now let $M\in {}_A\Mm^H$
If this property holds in the category ${}_A\Mm^H$, then we call $M$ an $H$-progenerator.

\begin{corollary}\colabel{5.8}
Assume that $A$ and $B$ are right faithfully flat $H$-Galois extensions. If $(A,B,M,N,\alpha,\beta)$
is a strict $H$-Morita context, then $M$ is an $H$-progenerator.
\end{corollary}

\begin{proof}
Let $(A^{{\rm co}H},B^{{\rm co}H},M_1,N_1,\alpha_1,\beta_1)$ be the corresponding strict
$\sq_H$-Morita context, as in \thref{5.7}. It follows from classical Morita theory that
$M_1$ is a left $A^{{\rm co}H}$-progenerator. The claim then follows from the equivalence
$A\ot_{A^{{\rm co}H}}-$ between the categories ${}_{A^{{\rm co}H}}\Mm$ and ${}_A\Mm^H$.
\end{proof}

\begin{theorem}\thlabel{5.9}
Assume that $A$ and $B$ are right faithfully flat $H$-Galois extensions, and let
$(A^{{\rm co}H},B^{{\rm co}H},M_1,N_1,\alpha_1,\beta_1)$ be a strict Morita context.
If $M_1$ has a left $A\sq_H B^{\rm op}$-module structure, then there is a unique
left $B\sq_H A^{\rm op}$-module structure on $N_1$ such that 
$(A^{{\rm co}H},B^{{\rm co}H},M_1,N_1,\alpha_1,\beta_1)$ is a strict
$\sq_H$-Morita context.
\end{theorem}

\begin{proof}
We know that $M=A\ot_{A^{{\rm co}H}}M_1\in {}_A\Mm_B^H$.
We have seen in \prref{3.7} that we have a morphism $\psi:\ B\to {}_A\END(M)$
in ${}_B\Mm_B^H$, and a left $B\sq_HB$-linear map
$\psi^{{\rm co}H}:\ B^{{\rm co}H}\to {}_A\END(M)^{{\rm co}H}\cong {}_{A^{{\rm co}H}}
\End(M^{{\rm co}H})$, see also \coref{3.4}. $\psi^{{\rm co}H}$ is an isomorphism, because
the Morita context is strict. Since $B$ is right faithfully flat $H$-Galois, it follows that
$\psi$ is an isomorphism in ${}_B\Mm_B^H$. Since $M_1$ is a progenerator as
a left $A^{{\rm co}H}$-module, $M$ is a progenerator as a left $A$-module.
Let $N={}_A\HOM(M,A)$. Then $N^{{\rm co}H}\cong {}_{A^{{\rm co}H}}\Hom(M_1,A^{{\rm co}H})$
as left $B\sq_HA^{\rm op}$-modules
(see \coref{3.4}); ${}_{A^{{\rm co}H}}\Hom(M_1,A^{{\rm co}H})$ and $N_1$
are canonically isomorphic as $(B^{{\rm co}H}\hbox{-}A^{{\rm co}H})$-bimodules, since the Morita context is strict. Using this isomorphism, the left $B\sq_HA^{\rm op}$-module structure
can be transported to $N_1$. The $H$-Morita context $(A,B,M,N)$ associated to $M$
is strict by \prref{5.4}. The corresponding $\sq_H$-Morita context from \thref{5.7}
is canonically isomorphic to $(A^{{\rm co}H},B^{{\rm co}H},M_1,N_1,\alpha_1,\beta_1)$.
This proves the claim.
\end{proof}

We end this Section with the following result.

\begin{theorem}\thlabel{5.10}
Let $A$ be a (right) faithfully flat Galois extension of $A^{{\rm co}H}$.
Assume that $M\in {}_A\Mm^H$ is a progenerator as a left $A$-module. Then $B={}_A\END(M)^{\rm op}$
is a (right) faithfully flat $H$-Galois extension of $B^{{\rm co}H}$ if and only if $M$ is an $H$-progenerator.
\end{theorem} 

\begin{proof}
The $H$-Morita context $(A,B,M,N={}_A\HOM(M,A),\alpha,\beta)$
from \exref{5.3} is strict by \prref{5.4}.\\
If $B$ is a faithfully flat $H$-Galois extension, then $M$ is an $H$-progenerator
by \coref{5.8}.\\
Conversely, let $M$ be an $H$-progenerator. $M\in {}_A\Mm_B^H$ (see \exref{5.3}),
hence $M_1=M^{{\rm co}H}\in {}_{A^{{\rm co}H}}\Mm_{B^{{\rm co}H}}$.
From the fact that the categories ${}_A\Mm^H$ and ${}_{A^{{\rm co}H}}\Mm$
are equivalent, it follows that $M_1$ is a left $A^{{\rm co}H}$-progenerator. From
\prref{3.3}, we know that $B^{{\rm co}H}\cong {}_{A^{{\rm co}H}}\End(M_1)^{\rm op}$
and that $N^{{\rm co}H}\cong {}_{A^{{\rm co}H}}\Hom(M_1, B^{{\rm co}H})$. The Morita context
$$(A^{{\rm co}H},B^{{\rm co}H}\cong {}_{A^{{\rm co}H}}\End(M_1)^{\rm op},
M_1,{}_{A^{{\rm co}H}}\Hom(M_1, B^{{\rm co}H}),\alpha_1,\beta_1)$$
associated to $M_1\in {}_{A^{{\rm co}H}}\Mm$ is strict, so $M_1\ot_{B^{{\rm co}H}}-
:\ {}_{B^{{\rm co}H}}\Mm\to  {}_{A^{{\rm co}H}}\Mm$ is a category equivalence.
$A\ot_{A^{{\rm co}H}}-:\ {}_{A^{{\rm co}H}}\Mm\to {}_A\Mm^H$ is an equivalence
since $A$ is a right faithfully flat $H$-Galois extension, and $M\ot_B-:\ {}_B\Mm^H\to
{}_A\Mm^H$ is also an equivalence (see \prref{5.2}). Using the fact that
$A\ot _{A^{{\rm co}H}}M_1\cong M$ ($A$ is a right  faithfully flat Galois extension), we
see easily that the following diagram of functors commutes;
$$\xymatrix{
{}_{B^{{\rm co}H}}\Mm\ar[rr]^{M_1\ot_{B^{{\rm co}H}}-}\ar[d]_{B\ot_{B^{{\rm co}H}}-}&&
{}_{B^{{\rm co}H}}\Mm\ar[d]^{A\ot_{A^{{\rm co}H}}-}\\
{}_B\Mm^H\ar[rr]^{M\ot_B-}&& {}_B\Mm^H}$$
Three of the four functors in the diagram are equivalences, hence the fourth one,
$B\ot_{B^{{\rm co}H}}-$ is also an equivalence (see the observations following
\coref{7.3}). $M_1$, $A$ and $M$ are right faithfully flat over $B^{{\rm co}H}$, $A^{{\rm co}H}$ and $B$ respectively, hence it follows that $B$ is right faithfully flat over $B^{{\rm co}H}$. Thus condition
(5) of \thref{1.1} is fulfilled, and it follows that $B$ is a right faithfully flat $H$-Galois extension.
\end{proof}

\section{Application to the Miyashita-Ulbrich action}\selabel{6}
Let $A$ be a right faithfully flat right $H$-Galois extension, and consider the map
$$\gamma_A=\can^{-1}\circ(\eta_A\ot H):\ H\to A\ot_{A^{{\rm co}H}}A.$$
Following \cite{Schauenburg}, we use the notation
$$\gamma_A(h)=\sum_i l_i(h)\ot_{A^{{\rm co}H}} r_i(h).$$
$\gamma_A(h)$ is then characterized by the property
$$\sum_i l_i(h)r_i(h)_{[0]}\ot r_i(h)_{[1]}=1\ot h.$$
The following properties are then easy to prove (see \cite[3.4]{Schneider1}):
for all $h,h'\in H$ and $a\in A$, we have
\begin{eqnarray}
&&\gamma_A(h)\in (A\ot_{A^{{\rm co}H}}A)^{A^{{\rm co}H}};\eqlabel{6.1.1}\\
&&\gamma_A(h_{(1)})\ot h_{(2)}=
\sum_i l_i(h)\ot_{A^{{\rm co}H}} r_i(h)_{[0]}\ot r_i(h)_{[1]};\eqlabel{6.1.2}\\
&&\gamma_A(h_{(2)})\ot S(h_{(1)})=
\sum_i l_i(h)_{[0]}\ot_{A^{{\rm co}H}} r_i(h)\ot l_i(h)_{[1]};\eqlabel{6.1.3}\\
&&\sum_i l_i(h)r_i(h)=\varepsilon(h)1_A;\eqlabel{6.1.4}\\
&&\sum_i a_{[0]}l_i(a_{[1]})\ot r_i(a_{[1]})=1\ot a;\eqlabel{6.1.5}\\
&&\gamma(hh')=\sum_{i,j} l_i(h')l_j(h)\ot_{A^{{\rm co}H}}
r_j(h)r_j(h').\eqlabel{6.1.6}
\end{eqnarray}
Combining \equref{6.1.2} and \equref{6.1.3}, we find
\begin{eqnarray}
&&\hspace*{-2cm}
\sum_i l_i(h)_{[0]}\ot_{A^{{\rm co}H}} r_i(h)_{[0]}\ot l_i(h)_{[1]}\ot r_i(h)_{[1]}\nonumber\\
&\equal{\equref{6.1.3}}&
\sum_i l_i(h_{(2)})\ot_{A^{{\rm co}H}} r_i(h_{(2)})_{[0]}
\ot S(h_{(1)}) \ot r_i(h_{(2)})_{[1]}\nonumber\\
&\equal{\equref{6.1.2}}& \sum_i l_i(h_{(2)})\ot_{A^{{\rm co}H}}
r_i(h_{(2)})\ot S(h_{(1)}) \ot h_{(3)}. \eqlabel{6.1.7}
\end{eqnarray}

Let $M$ be an $(A,A)$-bimodule. On $M^{A^{{\rm co}H}}$, we can define a right $H$-action
called the Miyashita-Ulbrich action. This was introduced in \cite{DT}, and we follow
here the description given in \cite{Schauenburg}. It is given by the formula
$$m\rightactr h= \sum_i l_i(h)mr_i(h).$$
It follows from \equref{6.1.1} and \equref{6.1.6} that we have a well-defined right
$H$-action. In particular, for $X,Y\in \Mm_A$, $\Hom(X,Y)\in {}_A\Mm_A$,
with left and right $A$-action given by
$$(a\cdot f\cdot a')(x)=f(xa)a'.$$
It is easy to see that $$\Hom(X,Y)^{A^{{\rm co}H}}=\Hom_{A^{{\rm
co}H}}(X,Y),$$ and the Miyashita-Ulbrich action is then given by
(see \cite[Cor. 3.5]{Schneider1})
$$(f \rightactr h)(x)=\sum_i f(xl_i(h))r_i(h).$$

\begin{lemma}\lelabel{6.1}
Let $A$ and $B$ be right faithfully flat right $H$-Galois extensions. For all $b\in B$,
we have that
$$x:=\gamma(S^{-1}(b_{[1]}))\ot b_{[0]}\in A\ot_{A^{{\rm co}H}}
(A\square_H B^{\rm op}).$$
\end{lemma}

\begin{proof} We have
\begin{eqnarray*}
&&\hspace*{-2cm}
\sum_i l_i(S^{-1}(b_{[1]})) \ot_{A^{{\rm co}H}} r_i(S^{-1}(b_{[1]}))_{[0]}
\ot r_i(S^{-1}(b_{[1]}))_{[1]}\ot b_{[0]}\\
&\equal{\equref{6.1.2}}& \gamma(S^{-1}(b_{[1]}))\ot S^{-1}(b_{[1]})\ot b_{[0]}.
\end{eqnarray*}
hence $x\in (A\ot_{A^{{\rm co}H}} A)\square_H B^{\rm op}\cong
A\ot_{A^{{\rm co}H}} (A\square_H B^{\rm op})$.
\end{proof}

Now we assume that $(A,B,M,N,\alpha,\beta)$ is a strict $H$-Morita context connecting
the right faithfully flat $H$-Galois extensions $A$ and $B$. For $X\in \Mm_A$, we
have the isomorphism
$$\diagram \varphi:\ X\ot_{A^{{\rm co}H}} M^{{\rm co}H}\cong X\ot_A A
\ot_{A^{{\rm co}H}} M^{{\rm co}H}\rrto^{\hspace{3cm}X\ot_A
\varepsilon_{1,M}} && X\ot_A M, \enddiagram$$ given by
$$\varphi(x\ot_{A^{{\rm co}H}}m)=x\ot_A m.$$ We have that $X\ot_A M\in \Mm_B$,
and its right $B$-action can be transported to $X\ot_{A^{{\rm
co}H}} M^{{\rm co}H}$. We compute this action in our next Lemma.

\begin{lemma}\lelabel{6.2}
The transported right $B$-action on $X\ot_{A^{{\rm co}H}} M^{{\rm co}H}$
is given by the formula
\begin{equation}\eqlabel{6.2.1}
(x\ot_{A^{{\rm co}H}} m)\cdot b= \sum_i xl_i(S^{-1}(b_{[1]}))\ot
_{A^{{\rm co}H}} (r_i(S^{-1}(b_{[1]}))\ot b_{[0]})m.
\end{equation}
\end{lemma}

\begin{proof}
Observe first that the action \equref{6.2.1} is well-defined, since
$M^{{\rm co}H}\in {}_{A\square B^{\rm op}}\Mm$, and by \leref{6.1}.
For the sake of simplicity, we introduce the following notation: for
$\sum_i a_i\ot b_i \in A\square_H B^{\rm op}$ and $m\in M^{{\rm co}H}$, we write
$$(\sum_i a_i\ot b_i)\cdot m=\sum_i a_imb_i.$$
We have to show that $\varphi$ is right $H$-linear. Indeed,
\begin{eqnarray*}
&&\hspace*{-2cm}
\varphi((x\ot_{A^{{\rm co}H}} m)\cdot b)=
\sum_i xl_i(S^{-1}(b_{[1]}))\ot _{A} r_i(S^{-1}(b_{[1]}))mb_{[0]}\\
&=&
\sum_i xl_i(S^{-1}(b_{[1]})) r_i(S^{-1}(b_{[1]}))\ot _{A} mb_{[0]}\\
&\equal{\equref{6.1.4}}&
\sum_i x\varepsilon(S^{-1}(b_{[1]}))\ot _{A} mb_{[0]}=x\ot_A mb.
\end{eqnarray*}
\end{proof}

Consider the setting of \thref{5.10}: $(A,B,M,N,\alpha,\beta)$ is a strict
$H$-Morita context connecting the right faithfully flat $H$-Galois extensions $A$ and $B$,
and $(A^{{\rm co}H},B^{{\rm co}H},M^{{\rm co}H},N^{{\rm co}H},\alpha_1,\beta_1)$
is the corresponding Morita context connecting $A^{{\rm co}H}$ and $B^{{\rm co}H}$.
For $X,Y\in \Mm_A$, we have an isomorphism
\begin{equation}\eqlabel{6.3.1}
\phi:\ \Hom_{A^{{\rm co}H}}(X,Y)\to \Hom_{B^{{\rm co}H}}(X\ot_{A^{{\rm co}H}}
M^{{\rm co}H}, Y\ot_{A^{{\rm co}H}} M^{{\rm co}H}),
\end{equation}
given by $\phi(f)= f\ot_{A^{{\rm co}H}} M^{{\rm co}H}$. It follows
from \leref{6.2} that $\Hom(X\ot_{A^{{\rm co}H}} M^{{\rm co}H},
Y\ot_{A^{{\rm co}H}} M^{{\rm co}H})$ is a $(B,B)$-bimodule, and we
can consider the Mi\-ya\-shi\-ta-Ulbrich action on $\Hom_{B^{{\rm
co}H}}(X\ot_{A^{{\rm co}H}} M^{{\rm co}H}, Y\ot_{A^{{\rm co}H}}
M^{{\rm co}H})$.

\begin{proposition}\prlabel{6.3}
With notation as above, the map $\phi$ from {\rm \equref{6.3.1}} preserves the
Miyashita-Ulbrich action.
\end{proposition}

\begin{proof}
We will use the notation
$$\gamma_B(h)=\sum_j k_j(h)\ot_{B^{{\rm co}H}} q_j(h)\in B\ot_{B^{{\rm co}H}} B.$$
We have to show that
$$\phi(f)\rightactr h=\phi(f\rightactr h),$$
for all right $A^{{\rm co}H}$-linear $f:\ X\to Y$ and $h\in H$. For $x\in X$ and
$m\in M^{{\rm co}H}$, we compute
\begin{eqnarray*}
&&\hspace*{-1cm}
(\phi(f)\rightactr h)(x\ot_{A^{{\rm co}H}} m)
= \sum_j \phi(f)\bigl( (x\ot_{A^{{\rm co}H}} m)k_j(h)\bigr) q_j(h)\\
&=&
\sum_{i,j}\Bigl[f\bigl[ xl_i(S^{-1}(k_j(h)_{[1]}))\bigr]
\ot_{A^{{\rm co}H}}
r_i(S^{-1}(k_j(h)_{[1]}))m k_j(h)_{[0]}\Bigr] q_j(h)\\
&=&
\sum_{i,j,p}f\bigl[ xl_i(S^{-1}(k_j(h)_{[1]}))\bigr] l_p(S^{-1}(q_j(h)_{[1]}))\\
&&\hspace*{1cm}\ot_{A^{{\rm co}H}}
r_p(S^{-1}(q_j(h)_{[1]}))r_i(S^{-1}(k_j(h)_{[1]}))m k_j(h)_{[0]}q_j(h)_{[0]}\\
&\equal{\equref{6.1.7}}&
\sum_{i,j,p}f\bigl[ xl_i(S^{-1}(S(h_{(1)})))\bigr] l_p(S^{-1}(h_{(3)}))\\
&&\hspace*{1cm}\ot_{A^{{\rm co}H}}
r_p(S^{-1}(h_{(3)}))r_i(S^{-1}(S(h_{(1)})))m k_j(h_{(2)})q_j(h_{(2)})\\
&\equal{\equref{6.1.4}}&
\sum_{i,p} f(xl_i(h_{(1)}))l_p(S^{-1}(h_{(2)}))
\ot_{A^{{\rm co}H}}r_p(S^{-1}(h_{(2)}))r_i(h_{(1)})m1_B\\
&\equal{(*)}&
\sum_{i,p} f(xl_i(h_{(1)}))l_p(S^{-1}(h_{(2)}))
r_p(S^{-1}(h_{(2)}))r_i(h_{(1)})\ot_{A^{{\rm co}H}} m\\
&\equal{\equref{6.1.4}}&
\sum_i f(xl_i(h))r_i(h) \ot_{A^{{\rm co}H}} m\\
&=& (f\rightactr h)(x) \ot_{A^{{\rm co}H}} m= (\phi(f\rightactr h))(x\ot_{A^{{\rm co}H}} m).
\end{eqnarray*}
The equality $(*)$ can be justified as follows. From \leref{6.1}, we deduce that, for
all $i$:
$$
\sum_i l_i(S^{-1}(k_j(h)_{[1]})) \ot_{A^{{\rm co}H}} 1_A
\ot_{A^{{\rm co}H}} r_i(S^{-1}(k_j(h)_{[1]}))\ot k_j(h)_{[0]}$$
and
$$
\sum_p 1_A  \ot_{A^{{\rm co}H}}l_p(S^{-1}(q_j(h)_{[1]}))
\ot_{A^{{\rm co}H}} r_p(S^{-1}(q_j(h)_{[1]}))\ot q_j(h)_{[0]}$$
lie in $(A\ot_{A^{{\rm co}H}}A)\ot_{A^{{\rm co}H}}
(A\square_HB^{\rm op})$. Consequently $(A\ot_{A^{{\rm
co}H}}A)\ot_{A^{{\rm co}H}} (A\square_HB^{\rm op})$ also contains
\begin{eqnarray*}
&&\hspace*{-1cm}
\sum_{i,p} l_i(S^{-1}(k_j(h)_{[1]})) \ot_{A^{{\rm co}H}}
l_p(S^{-1}(q_j(h)_{[1]}))\\
&&r_p(S^{-1}(q_j(h)_{[1]}))r_i(S^{-1}(k_j(h)_{[1]}))\ot k_j(h)_{[0]}q_j(h)_{[0]}
\end{eqnarray*}
\begin{eqnarray*}
&\equal{\eqref{eq:6.1.4},\eqref{eq:6.1.7}}&
\sum_{i,p} l_i(h_{(1)}) \ot_{A^{{\rm co}H}} l_p(S^{-1}(h_{(2)}))\ot_{A^{{\rm co}H}}
r_p(S^{-1}(h_{(2)}))r_i(h_{(1)}) \ot 1_B\\
&=:& Z\ot 1_B.
\end{eqnarray*}
This means that
$$(A\ot_{A^{{\rm co}H}} A\ot_{A^{{\rm co}H}}\rho_A)(Z)\ot 1_B=
Z\ot 1_H\ot 1_B,$$ hence
$$Z\in (A\ot_{A^{{\rm co}H}} A\ot_{A^{{\rm co}H}} A)^{{\rm co}H}\cong
A\ot_{A^{{\rm co}H}} A\ot_{A^{{\rm co}H}}A^{{\rm co}H},
$$
since $A/A^{{\rm co}H}$ is faihfully flat.
\end{proof}

\section{Hopf subalgebras}\selabel{7}
Throughout this Section, $H$ is a Hopf algebra with bijective antipode over a field $k$,
and $K$ is a Hopf subalgebra of $H$. We assume that the antipode of $K$ is bijective,
and that $H$ is (right) faithfully flat as a left $K$-module. Let $K^+=\Ker(\varepsilon_K)$.
It is well-known, and easy to prove (see \cite[Sec. 1]{Ulbrich}) that
$$\ol{H}=H/HK^+\cong H\ot_K k$$
is a left $H$-module coalgebra, with operations
$$h\cdot \ol{l}=\ol{hl},~~\Delta_{\ol{H}}(\ol{h})=\ol{h}_{(1)}\ot \ol{h}_{(2)},~~
\varepsilon_{\ol{H}}(\ol{h})=\varepsilon(h).$$
The class in $\ol{H}$ represented by $h\in H$ is denoted by $\ol{h}$.
$\ol{1}$ is a grouplike element of $\ol{H}$, and we consider coinvariants with
respect to this element. A right $H$-comodule $M$ is also a right $\ol{H}$-comodule,
by corestriction of coscalars:
$$\rho_{\ol{H}}(m)=m_{[0]}\ot \ol{m}_{[1]}.$$
The $\ol{H}$-coinvariants of $M\in \Mm^H$ are then
\begin{eqnarray*}
M^{{\rm co}\ol{H}}&=&\{m\in M~|~m_{[0]}\ot \ol{m}_{[1]}=m\ot \ol{1}\}\\
&=& \{m\in M~|~\rho(m)\in M\ot K\}\cong M\square_H K.
\end{eqnarray*}
If $A$ is a right $H$-comodule algebra, then $A^{{\rm co}\ol{H}}$
is a right $K$-comodule algebra, and $(A^{{\rm co}\ol{H}})^{{\rm co}K}
=A^{{\rm co}H}$.

\begin{proposition}\prlabel{7.1} {\bf (\cite[Remark 1.8]{Schneider1})}
Let $H$, $K$ and $A$ be as above, and assume that $A$ is a faithfully
flat $H$-Galois extension. Then $A$ is right faithfully flat as a right $A^{{\rm co}\ol{H}}$-module,
and
$$\can:\ A\ot_{A^{{\rm co}\ol{H}}} A\to A\ot \ol{H},~~
\can (a\ot b)=ab_{[0]}\ot \ol{b}_{[1]}$$
is bijective. The functors $(A\ot_{A^{{\rm co}\ol{H}}}-,(-)^{{\rm co}\ol{H}})$
form a pair of inverse equivalences between the categories
${}_{A^{{\rm co}\ol{H}}}\Mm$ and ${}_A\Mm(H)^{\ol{H}}$.
\end{proposition}

We also have an adjoint pair of functors $(F_4=A\ot_{A^{{\rm co}\ol{H}}}-,
G_4=(-)^{{\rm co}\ol{H}}\cong -\square_H K)$ between the categories
 ${}_{A^{{\rm co}\ol{H}}}\Mm^K$ and ${}_A\Mm^H$. This can be seen directly,
 but it is also a consequence of a more general result: we apply
 \cite[Theorem 1.3]{CaenepeelR} to the inclusion morphism between the
 Doi-Hopf data $(K, A^{{\rm co}\ol{H}}, K)$ and $(H,A,H)$.\\
 Take $N\in {}_{A^{{\rm co}\ol{H}}}\Mm^K$. Forgetting the $K$-coaction, we find
 that $N\in  {}_{A^{{\rm co}\ol{H}}}\Mm$. Then it is easy to see that the counit map
 $\eta_N:\ N\to (A\ot_{A^{{\rm co}\ol{H}}}N^{{\rm co}\ol{H}}$ is a morphism in
 ${}_{A^{{\rm co}\ol{H}}}\Mm^K$, and coincides with the counit map from the
 adjunction $(F_4,G_4)$. Since $\eta_N$ is an isomorphism, the unit maps
 of the adjunction $(F_4,G_4)$ are isomorphisms. In the same way, we can conclude
 that the counit maps are isomorphisms, and we conclude
 
 \begin{corollary}\colabel{7.3}
 Let $H$, $K$ and $A$ be as above, and assume that $A$ is a faithfully
flat $H$-Galois extension.
 Then the adjoint pair of functors
 $(F_4=A\ot_{A^{{\rm co}\ol{H}}}-,
G_4=(-)^{{\rm co}\ol{H}}\cong -\square_H K)$ establishes a pair of inverse
equivalences between the categories ${}_{A^{{\rm co}\ol{H}}}\Mm^K$ and ${}_A\Mm^H$.
\end{corollary}

Before stating our next corollary, we recall some elementary facts from category theory.
If $(F_1,G_1)$ and $(F_2,G_2)$ are pairs of adjoint functors, respectively
between categories $\Cc$ and $\Dd$, and between $\Dd$ and $\Ee$,
then $(F=F_2\circ F_1,G=G_1\circ G_2)$ is a pair of adjoint functors between
$\Cc$ and $\Ee$. If two of these three pairs are inverse equivalences, then the
third one is also a pair of inverse equivalences. This follows from the following
formulas, which give the relations between the units and counits of the three
adjunctions: for all $C\in \Cc$ and $D\in \Dd$, we have
$$\eta_C=G_1(\eta_{2,F_1(C)})\circ \eta_{1,C}~~;~~
\varepsilon_E= \eta_{2,E}\circ F_2(\varepsilon_{1,G_2(E)}).$$

This can be applied to the following situation. Assume that we are in the setting of
\prref{7.1} and \coref{7.3}. We have adjunctions
\begin{itemize}
\item $(F_1=A\ot_{A^{{\rm co}H}}-,\ G_1=(-)^{{\rm co}H})$ between
${}_{A^{{\rm co}H}}\Mm$ and ${}_A\Mm^H$; \item $(F_3=A^{{\rm
co}\ol{H}}\ot_{A^{{\rm co}H}}-,\ G_3=(-)^{{\rm co}K})$ between
${}_{A^{{\rm co}H}}\Mm$ and ${}_{A^{{\rm co}\ol{H}}}\Mm^K$;
\item $(F_4=A\ot_{A^{{\rm co}\ol{H}}}-,
G_4=(-)^{{\rm co}\ol{H}}\cong -\square_H K)$ between the categories ${}_{A^{{\rm co}\ol{H}}}\Mm^K$ and ${}_A\Mm^H$.
\end{itemize}
It is clear that $F_1=F_4\circ F_3$ and $G_1=G_3\circ G_4$.
$(F_1,G_1)$ and $(F_4,G_4)$ are pairs of inverse equivalences, by \thref{1.1}
and \coref{7.3}. Hence $(F_3,G_3)$ is also a pair of inverse equivalences, and using
\thref{1.1}, we obtain the following result.

\begin{corollary}\colabel{7.2}
Let $H$, $K$ and $A$ be as above, and assume that $A$ is a faithfully
flat $H$-Galois extension. Then $A^{{\rm co}\ol{H}}$ is a right faithfully flat
$K$-Galois extension.
\end{corollary}

\begin{theorem}\thlabel{7.4}
Let $H$ and $K$ be as before: $K\subset H$ are Hopf algebras with invertible
antipode over a field $k$, and $H$ is faithfully flat as a left $K$-module.
Let $A$ and $B$ be (right) faithfully flat right $H$-Galois extensions, connected by
a strict $H$-Morita context $(A, B, M,N, \alpha,\beta)$.
\begin{enumerate}
\item $A^{{\rm co}\ol{H}}$ and $B^{{\rm co}\ol{H}}$ are connected
by a strict $K$-Morita context, with connecting modules $M^{{\rm
co}\ol{H}}$ and $N^{{\rm co}\ol{H}}$; \item we have a pair of
inverse equivalences $(M\ot_B-, N\ot_A-)$ between the categories
${}_B\Mm(H)^{\ol{H}}$ and ${}_A\Mm(H)^{\ol{H}}$; \item the
following diagram of categories and functors  commutes to within
natural equivalences of functors:
$$\diagram {}_A\mathcal{M}(H)^{{\bar H}}\rrto<3pt>^{N\otimes_A-}\ddto<3pt>^{(-){{^{\operatorname{co}{\bar H}}}}}
&& _B\mathcal{M}(H)^{{\bar H}}\llto<3pt>^{M\otimes_B-} \ddto<3pt>^{(-){{^{\operatorname{co}{\bar H}}}}} \\ \\
_{A{{^{\operatorname{co}{\bar
H}}}}}\mathcal{M}\rrto<3pt>^{N{{^{\operatorname{co}{\bar
H}}}}\otimes_{A{{^{\operatorname{co}{\bar
H}}}}}-}\uuto<3pt>^{A\otimes_{A{{^{\operatorname{co}{\bar H}}}}}-}
&& _{B{{^{\operatorname{co}{\bar
H}}}}}\mathcal{M}\llto<3pt>^{M{{^{\operatorname{co}{\bar
H}}}}\otimes_{B{{^{\operatorname{co}{\bar
H}}}}}-}\uuto<3pt>^{B\otimes_{B{{^{\operatorname{co}{\bar
H}}}}}-}\enddiagram$$
\end{enumerate}
\end{theorem}

\begin{proof} 
1) We have the following commutative diagram of inclusions
$$\diagram
A^{{\rm co}\ol{H}}\square_K (B^{\rm op})^{{\rm co}\ol{H}}\rrto^{\subset}
\dto_{\subset}&&A^{{\rm co}\ol{H}}\ot (B^{\rm op})^{{\rm co}\ol{H}}
\dto^{\subset}\\
A\square_H B^{\rm op}\rrto^{\subset} &&A\ot B^{\rm op}
\enddiagram$$
By \thref{5.7}, we have
a strict $\sq_H$-Morita context $(A^{{\rm co}H},B^{{\rm co}H},M^{{\rm co}H},$ $N^{{\rm co}H},
\alpha_1,\beta_1)$. By restriction of scalars, A is a left $A^{{\rm co}\ol{H}}\square_K (B^{\rm op})^{{\rm co}\ol{H}}$-module. Then we can apply Theorems \ref{th:5.7} and \ref{th:5.9},
with $H$ replaced by $K$, and taking into account that $A^{{\rm co}\ol{H}}$
and $B^{{\rm co}\ol{H}}$ are right faithfully flat $K$-Galois extensions, by \coref{7.2}.
We find that $A^{{\rm co}\ol{H}}$ and $B^{{\rm co}\ol{H}}$ are connected by
a strict $K$-Morita context. The first connecting module is
\begin{eqnarray*}
&&\hspace*{-2cm}
A^{{\rm co}\ol{H}}\ot_{A^{{\rm co}H}} M^{{\rm co}H}
=F_3G_1(M)\cong G_4F_4F_3G_1(M)\\
&\cong& G_4F_1G_1(M)\cong G_4(M)=M^{{\rm co}\ol{H}}.
\end{eqnarray*}
In a similar way, we find that the second connecting module is $N^{{\rm co}\ol{H}}$.\\
2) The proof is an easy adaption of the proof of \prref{5.2}.\\
3) $B^{{\rm co}\ol{H}}$ is a right $K$-comodule algebra, and, by corestriction of coscalars,
a right $H$-comodule algebra, so we can consider the categories
${}_{A^{{\rm co}\ol{H}}}\Mm^K_{B^{{\rm co}\ol{H}}}$ and ${}_A\Mm^H_{B^{{\rm co}\ol{H}}}$.
It is then easy to see that the inverse equivalent functors of \coref{7.3}
also define a pair of inverse equivalences between these two categories of relative Hopf
bimodules. Now $M\in {}_A\Mm^H_{B^{{\rm co}\ol{H}}}$, so
$M\cong A\ot_{A^{{\rm co}\ol{H}}} M^{{\rm co}\ol{H}}$ as right $B^{{\rm co}\ol{H}}$-modules.
It follows that we have, for all $P\in {}_{B^{{\rm co}\ol{H}}}\Mm$,
$$A\ot_{A^{{\rm co}\ol{H}}} M^{{\rm co}\ol{H}}\ot_{B^{{\rm co}\ol{H}}}  P
\cong M\ot_{B^{{\rm co}\ol{H}}}  P\cong M\ot_BB\ot_{B^{{\rm co}\ol{H}}}  P.$$
In a similar way, we can show that
$$B\ot_{B^{{\rm co}\ol{H}}} N^{{\rm co}\ol{H}}\ot_{A^{{\rm co}\ol{H}}}  Q\cong
N\ot_AA\ot_{A^{{\rm co}\ol{H}}} Q,$$
for all $Q {}_{A^{{\rm co}\ol{H}}}\Mm$. Finally, take $U\in {}_A\Mm(H)^{\ol{H}}$. Then
\begin{eqnarray*}
&&\hspace*{-2cm}
(N\ot_A U)^{{\rm co}\ol{H}}\cong (N\ot_AA\ot_{A^{{\rm co}\ol{H}}} U^{{\rm co}\ol{H}})^{{\rm co}\ol{H}}\\
&\cong& (B\ot_{B^{{\rm co}\ol{H}}}N^{{\rm co}\ol{H}}\ot_{A^{{\rm co}\ol{H}}} U^{{\rm co}\ol{H}})^{{\rm co}\ol{H}}\cong N^{{\rm co}\ol{H}}\ot_{A^{{\rm co}\ol{H}}} U^{{\rm co}\ol{H}},
\end{eqnarray*}
and, in a similar way, for $V\in {}_B\Mm(H)^{\ol{H}}$,
$$(M\ot_B V)^{{\rm co}\ol{H}}\cong
M^{{\rm co}\ol{H}}\ot_{B^{{\rm co}\ol{H}}} V^{{\rm co}\ol{H}}.$$
\end{proof}

Finally recall that if the algebras $A$ and $B$ are Morita
equivalent, then there is a Morita equivalence between $A\otimes
A{{^{\operatorname{op}}}}$ and $B\otimes
B{{^{\operatorname{op}}}}$ sending $A$ to $B$. In particular, this
implies that the centers of $A$ and $B$ are isomorphic. In our
context this generalizes as follows.

\begin{corollary}\colabel{7.5}
Assume that $(A,B,M,N,\alpha,\beta)$ is a strict $H$-Morita context.
\begin{enumerate}
\item Let  $K$ and $L$ be  Hopf subalgebras of $H$ with
bijective antipodes, and assume that $H\otimes H$ is faithfully
flat as a right $K\otimes L$-module. Then the categories
$_{A{{^{\operatorname{co}H/HK^+}}}}\mathcal{M}_{A{{^{\operatorname{co}H/HL^+}}}}$
and
$_{B{{^{\operatorname{co}H/HK^+}}}}\mathcal{M}_{B{{^{\operatorname{co}H/HL^+}}}}$
are equivalent.
\item There is an isomorphism
$$C_A(A{{^{\operatorname{co}H}}})\cong
C_B(B{{^{\operatorname{co}H}}})$$ of left $H$-module right
$H$-comodule algebras, where  $C_A(A{{^{\operatorname{co}H}}})$ denotes the centralizer in $A$ of $A{{^{\operatorname{co}H}}}$.
\end{enumerate}
\end{corollary}

\begin{proof}
1) The objects $M\otimes N\in \ _{A\otimes A{{^{\operatorname{op}}}}}\mathcal{M}^{H\otimes H}_{B\otimes
B{{^{\operatorname{op}}}}}$ and $N\otimes M\in \ _{B\otimes
B{{^{\operatorname{op}}}}}\mathcal{M}^{H\otimes H}_{A\otimes
A{{^{\operatorname{op}}}}}$ induce a Morita equivalence between
$A\otimes A{{^{\operatorname{op}}}}$ and $B\otimes
B{{^{\operatorname{op}}}}$. Now the assertion follows from
\thref{7.4}, where we replace $H$ by $H\otimes H$,
$K$ by $K\otimes L$, $A$ by $A\otimes A{{^{\operatorname{op}}}}$
and $B$ by $B\otimes B{{^{\operatorname{op}}}}$.\\
2) Note that $$C_A(A{{^{\operatorname{co}H}}})\cong
\End_{A{{^{\operatorname{co}H}}}\otimes
A{{^{\operatorname{op}}}}}(A)$$ as $H$-module $H$-comodule
algebras. Since  under the equivalence of 1) (where we take $K=k$
and $L=H$), $A$ corresponds to $B$, the statement follows from
\prref{6.3}.
\end{proof}

\section{$H$-colinear equivalences}\selabel{8}
Let $H$ be a projective Hopf algebra, and $A$ a right $H$-comodule algebra.
Let ${}_A\ul{\Mm}^H$ be the category with relative Hopf modules as modules;
the set of morphisms between two objects $M$ and $N$ is ${}_A\HOM(M,N)$.
${}_A\ul{\Mm}^H$ is a right $H$-colinear category in the following sense:
${}_A\HOM(M,N)$ is a right $H$-comodule (see \prref{3.2x}); the map
$$\varphi:\ M\ot {}_A\HOM(M,N)\to N,~~\varphi(m\ot f)=f(m)$$
is right $H$-colinear (take $B=C=k$ in \prref{3.5}); if $N$ is a third object in
${}_A\ul{\Mm}^H$, then the composition
$$\psi:\ {}_A\HOM(L,M)\ot {}_A\HOM(M,N)\to {}_A\HOM(L,N),~~
\psi(f\ot g)=g\circ f$$
is right $H$-colinear. The following result is then obvious.

\begin{proposition}\prlabel{8.1}
Let $H$ be a projective Hopf algebra.
Let $(A,B,M,N,\alpha,\beta)$ be a strict $H$-Morita context connecting the
right $H$-comodule algebras $A$ and $B$. Then the functors $M\ot_B-$
and $N\ot_A-$ induce a pair of inverse right $H$-colinear equivalences
between  ${}_B\ul{\Mm}^H$ and ${}_A\ul{\Mm}^H$.
\end{proposition}

The functors $F=M\ot_B-$ and $G=N\ot_A-$ are right $H$-colinear in the following
sense: for $V,W\in {}_B\Mm^H$, the map
$$F:\ {}_B\HOM(V,W)\to {}_A\HOM(M\ot_B V,M\ot_BW),~~F(f)=M\ot_B f$$
is right $H$-colinear.\\
In this Section, we investigate when the converse of \prref{8.1} holds: suppose that we have
a pair of inverse right $H$-colinear equivalences
between  ${}_B\ul{\Mm}^H$ and ${}_A\ul{\Mm}^H$. Is this equivalence induced by
a strict $H$-Morita context? 
To this end, we will give an $H$-colinear version of the Eilenberg-Watts Theorem.

\begin{proposition}\prlabel{8.4}
Let $A$ and $B$ be $H$-comodule algebras, and $T:\ {}_{A}\ul{\Mm}^H\to  {}_{B}\ul{\Mm}^H$
an $H$-colinear functor. Then $N=T(A)\in {}_A\Mm_B^H$,
and we have a natural transformation $\psi:\ F=N\ot_A-\to T$, such that
$\psi_A:\ N\ot_AA\to T(A)=N$ is the natural isomorphism.
\end{proposition}

\begin{proof}
In the sequel, $V$ and $W$ will be objects in  ${}_{A}\ul{\Mm}^H$.
The fact that $T$ is right $H$-colinear means that
\begin{equation}\eqlabel{8.4.a}
T(f_{[0]})\ot f_{[1]}=\rho(T(f)),
\end{equation}
for $f\in {}A\HOM(V,W)$. We claim that the map
$$\varphi_V:\ V\to {}_A\HOM(A,V),~~\varphi_V(v)(a)=av$$
Is well-defined and right $H$-colinear. To this end, it suffices to show that
\begin{equation}\eqlabel{8.4.b}
\varphi(v)_{[0]}\ot \varphi(v)_{[1]}=\varphi(v_{[0]})\ot v_{[1]},
\end{equation}
for all $v\in V$. For all $a\in A$, we have
\begin{eqnarray*}
&&\hspace*{-2cm}
\bigl(\varphi_V(v)(a_{[0]})\bigr)_{[0]}\ot S^{-1}(a_{[1]})\bigl(\varphi_V(v)(a_{[0]})\bigr)_{[1]}\\
&=& a_{[0]}v_{[0]}\ot S^{-1}(a_{[2]})a_{[1]}v_{[1]}\\
&=& av_{[0]}\ot v_{[1]}=\varphi_V(v_{[0]})(a)\ot v_{[1]},
\end{eqnarray*}
and \equref{8.4.b} follows using \equref{3.2.1}.
$\varphi_V$ satisfies the following property:
\begin{equation}\eqlabel{8.4.1}
\varphi_V(av)=\varphi_V(v)\circ \varphi_A(a),
\end{equation}
for all $a\in A$ and $v\in V$. Indeed,
$$\varphi_V(av)(c)=cav=\varphi_V(v)(ca)=(\varphi_V(v)\circ \varphi_A(a))(c).$$
On $N=T(A)\in {}_{B}\Mm^H$, we define a right $A$-action as follows:
$$na=T(\varphi_A(a))(n),$$
for all $a\in A$ and $n\in N$. This makes $N$ an object of ${}_B\Mm_A^H$,
since
\begin{eqnarray*}
&&\hspace*{-2cm}
n(ac)= T(\varphi_A(ac))(n)\equal{\equref{8.4.1}} T(\varphi_A(a)\circ \varphi_A(c))(n)\\
&=& \bigl( T(\varphi_A(a))\circ T(\varphi_A(c))\bigr)(n)=(na)c;\\
&&\hspace*{-2cm}
(bn)a=T(\varphi_A(a))(bn)=bT(\varphi_A(a))(n)=b(na);\\
&&\hspace*{-2cm}
n_{[0]}a_{[0]}\ot n_{[1]}a_{[1]}= T(\varphi(a_{[0]}))(n_{[0]}) \ot n_{[1]}a_{[1]}\\
&\equal{\equref{8.4.b}}& T(\varphi(a)_{[0]})(n_{[0]}) \ot n_{[1]}\varphi(a)_{[1]}\\
&\equal{\equref{8.4.a}}& T(\varphi(a))_{[0]}(n_{[0]}) \ot n_{[1]}T(\varphi(a))_{[1]}\\
&\equal{\equref{3.2.1}}& T(\varphi(a))(n_{[0]})_{[0]}\ot n_{[2]}S^{-1}(n_{[1]})
T(\varphi(a))(n_{[0]})_{[1]}\\
&=& T(\varphi(a))(n)_{[0]}\ot T(\varphi(a))(n)_{[1]}=\rho(na),
\end{eqnarray*}
for all $a,c\in A$, $b\in B$ and $n\in N$.\\

For every $v\in V$, $\varphi_V(v):\ A\to V$ is left $A$-linear, hence $T(\varphi_V(v)):\
T(A)=N\to T(V)$ is left $B$-linear. By (\ref{eq:8.4.a},\ref{eq:8.4.b}), we also have that
\begin{equation}\eqlabel{8.4.2}
T(\varphi_V(v_{[0]}))\ot v_{[1]}=\rho(T(\varphi_V(v))).
\end{equation}
Now we define
$$\psi_V:\ N\ot_A V\to T(V),~~\psi_V(n\ot_A v)= T(\varphi_V(v))(n).$$
$\psi_V$ is well-defined since
\begin{eqnarray*}
&&\hspace*{-2cm}
\psi_V(n\ot_A av)=T(\varphi_V(av))(n)\equal{\equref{8.4.1}}
(T(\varphi_V(v))\circ T(\varphi_V(a)))(n)\\
&=& T(\varphi_V(v))(na)=
\psi_V(n\ot_A v).
\end{eqnarray*}
$\psi_V$ is right $H$-colinear, since
\begin{eqnarray*}
&&\hspace*{-2cm}
\psi_V(n_{[0]}\ot_A v_{[0]})\ot n_{[1]}v_{[1]}= T(\varphi_V(v_{[0]}))(n_{[0]})\ot n_{[1]}v_{[1]}\\
&\equal{\equref{8.4.2}}& T(\varphi_V(v))_{[0]}(n_{[0]})\ot n_{[1]}T(\varphi_V(v))_{[1]}\\
&\equal{\equref{3.2.1}}&\bigl(T(\varphi_V(v))(n_{[0]})\bigr)_{[0]}
\ot n_{[2]}S^{-1}(n_{[1]}) \bigl(T(\varphi_V(v))(n_{[0]})\bigr)_{[1]}\\
&=& \rho\bigl(T(\varphi_V(v))(n)\bigr)=\rho(\psi_V(n\ot_A v)).
\end{eqnarray*}
$\psi_V$ is left $B$-linear, since
$$\psi_V(bn\ot_A v)=T(\varphi_V(v))(bm)=b(T(\varphi_V(v))(m))=b\psi_V(n\ot_A v).$$
In order to show that $\psi$ is a natural transformation, we first observe the following
property. For $f:\ V\to W$ in ${}_{A}\Mm^H$, $v\in V$ and $a\in A$, we have
$$\varphi_W(f(v))(a)=af(v)=f(av)=f(\varphi_V(v)(a),$$
so $\varphi_W(f(v))=f\circ \varphi_V(v)$. We can now show that the diagram
$$\diagram
N\ot_AV\rrto^{\psi_V}\dto_{N\ot_A f}&&T(V)\dto^{T(f)}\\
N\ot_AW\rrto^{\psi_W}&&T(W)
\enddiagram$$
commutes:
\begin{eqnarray*}
&&\hspace*{-2cm}
(T(f)\circ \psi_V)(n\ot_A v)=(T(f)\circ T(\varphi_V(v)))(n)\\
&=& T(f\circ\varphi_V(v))(n)=T(\varphi_W(f(v)))(n)=\psi_W(n\ot_A f(v))).
\end{eqnarray*}
It follows that $\psi$ is a natural transformation. Finally, it is easy to compute that the map
$\psi_A:\ N\ot_A \to T(A)=A$ is given by
$\psi_A(n\ot a)=T(\varphi_A(a))(n)=na$,
as needed.
\end{proof}

We are now ready to prove the following generalization of the Eilenberg-Watts
Theorem (cf. \cite[II.2.3]{Bass}).

\begin{proposition}\prlabel{8.5}
With notation and assumptions as in \prref{8.4}, assume that $A$ is a generator
of ${}_{A}\Mm^H$, and that $T$, viewed as a functor ${}_{A}\Mm^H\to {}_{B}\Mm^H$,
preserves cokernels and arbitrary coproducts. Then the natural transformation
$\psi:\ F=N\ot_A-\to T$ from \prref{8.4} is a natural isomorphism.
\end{proposition}

\begin{proof}
Let $I$ be an index set, and $A^{(I)}$ the coproduct of copies of $A$ indexed by $I$.
For $i\in I$, let $r_i:\ A\to A^{(I)}$ be the natural inclusion. Since $\psi$ is a natural
transformation, we have a commutative diagram
$$\diagram
F(A)\rrto^{\psi_A}\dto_{F(r_i)}&&T(A)\dto^{T(r_i)}\\
F(A^{(I)})\rrto^{\psi_{A^{(I)}}}&&T(A^{(I)})
\enddiagram$$
Let $n_i:\ T(A)\to T(A)^{(I)}$ be the natural inclusion. Then the diagram
$$\diagram
F(A)^{(I)}\rrto^{\oplus_{i\in I}(n_i\circ\psi_A)}\dto_{\oplus_{i\in I}F(r_i)}&&T(A)^{(I)}\dto^{\oplus_{i\in I}T(r_i)}\\
F(A^{(I)})\rrto^{\psi_{A^{(I)}}}&&T(A^{(I)})
\enddiagram$$
also commutes. The vertical maps in the diagram are isomorphisms, since $F$ and $T$
commute with direct sums. We have seen in \prref{8.4} that the top horizontal map is an isomorphim,
so it follows that $\psi_{A^{(I)}}$ is an isomorphism.\\
Now take an arbitrary $V\in {}_{A}\Mm^H$. Since $A$ is a generator of ${}_{A}\Mm^H$,
we have an exact sequence
$$A^{(J)}\mapright{\pi} A^{(I)}\mapright{\varphi} V\to 0$$
in ${}_{A}\Mm^H$. Since $\psi$ is a natural transformation, and $F$ and $G$ preserve
cokernels, we have the following commutative diagram with exact rows in ${}_{B}\Mm^H$:
$$\diagram
F(A^{(J)})\dto^{\psi_{A^{(J)}}}\rrto^{F(\pi)}&&
F(A^{(I)})\dto^{\psi_{A^{(I)}}}\rrto^{F(\varphi)}&&
F(V)\dto^{\psi_V}\rrto^{}&&0\\
T(A^{(J)})\rrto^{T(\pi)}&&T(A^{(I)})\rrto^{T(\varphi)}&&
T(V)\rrto^{}&&0
\enddiagram$$
We know from above that $\psi_{A^{(J)}}$ and $\psi_{A^{(I)}}$ are isomorphisms, and it
follows from the 5 lemma that $\psi_V$ is also an isomorphism.
\end{proof}

\begin{theorem}\thlabel{8.6}
Let $A$ and $B$ be $H$-module algebras, and suppose that they generate the categories
${}_{A}\Mm^H$ and ${}_{B}\Mm^H$. If $(T,U)$ is a pair of $H$-linear inverse equivalences between
the categories ${}_{A}\ul{\Mm}^H$ and ${}_{B}\ul{\Mm}^H$, then there exists a strict $H$-Morita 
context $(A,B,M,N,\alpha,\beta)$ such that $T\cong N\ot_A-$ and $U\cong M\ot_B-$.
\end{theorem}

\begin{proof}
Since $(T,U)$ is also a pair of inverse equivalences between ${}_{A}\Mm^H$ and ${}_{B}\Mm^H$,
$T$ and $U$ preserve coproducts and cokernels. Applying \prref{8.5}, we find
$M\in {}_A\Mm_B^H$ and $N\in {}_B\Mm_A^H$
such that $T\cong N\ot_A-$ and $U\cong M\ot_B-$.\\
$(T,U)$ is a pair of adjoint functors, and the unit $\eta$ and the counit $\varepsilon$ are
natural isomorphisms. We define $\alpha=\eta^{-1}_A:\ M\ot_B N\to A$. Then $\alpha\in
{}_{A}\Mm^H$. Let us show that $\alpha$ is also right $A$-linear. For every $c\in A$,
the map $f_c:\ A\to A$, $f_c(a)=ac$ is left $A$-linear. Since $\eta$ is a natural transformation,
the diagram
$$\xymatrix{
A\ar[rr]^-{\eta_A}\ar[d]_{f_c}&&M\ot_BN\ot_AA\ar[d]^{M\ot_BN\ot_Af_c}\\
A\ar[rr]^-{\eta_A}&&M\ot_BN\ot_AA
}$$
commutes. Evaluating the diagram at $1_A$, we find that $\eta_A(ac)=\eta_A(a)c$.\\
We define $\beta=\varepsilon_B:\ N\ot_A M\to B$. Applying the above argument to the
adjunction $(U,T)$ with unit $\varepsilon^{-1}$ and counit $\eta^{-1}$, we find that
$\varepsilon_B$ is right $B$-linear.\\
Take $W\in {}_{B}\Mm^H$. For every $w\in W$, we consider the left $B$-linear map
$g_w:\ B\to W$, $g_w(b)=bw$. Since $\varepsilon$ is a natural transformation, the
diagram
$$\xymatrix{
N\ot_A M\ot_B B\ar[rr]^-{\varepsilon_B}\ar[d]_{N\ot_A M\ot_B g_w}&&B\ar[d]^{g_w}\\
N\ot_A M\ot_B W\ar[rr]^-{\varepsilon_W}&&W
}$$
commutes. Evaluating the diagram at $n\ot_A m\ot_B 1$, we see that
$\varepsilon_W=\varepsilon_B\ot_B W$.\\
From the properties of adjoint functors, we know that
$\varepsilon_{T(V)}\circ T(\eta_V)=T(V)$, for all $V\in {}_{A}\Mm^H$ . Taking
$V=A$ in this formula, we see that the diagram
$$\xymatrix{
N\ot_A\ar[rr]^-{N\ot_A\eta_A}\ar[d]_{\cong}&& N\ot_A M\ot_B N\ar[d]^{\varepsilon_N=\varepsilon_B\ot_BN}\\
N\ar[rr]^-{\cong}&&B\ot_B N
}$$
commutes. This diagram is one of the two diagrams in the definition of a Morita context.
The commutativity of the other diagram follows in a similar way.
\end{proof}

\begin{corollary}\colabel{8.7}
Let $H$ be a projective Hopf algebra, and assume that the right $H$-comodule
algebras $A$ and $B$ are $H$-Galois extensions of $A^{{\rm co}H}$ and $B^{{\rm co}H}$,
respectively. If $(T,U)$ is a pair of $H$-colinear inverse equivalences between
the categories ${}_{A}\ul{\Mm}^H$ and ${}_{B}\ul{\Mm}^H$, then there exists a strict $H$-Morita context $(A,B,M,N,\alpha,\beta)$ such that $T\cong N\ot_A-$ and $U\cong M\ot_B-$.
\end{corollary}

\begin{proof}
It is well-known that $A^{{\rm co}H}$ is a generator of ${}_{A^{{\rm co}H}}\Mm$; since $(F_1,G_1)$
is a pair of inverse equivalences (see \thref{1.1}), $F_1(A^{{\rm co}H})=A$ is a generator of
${}_A\Mm^H$. In a similar way, $B$ is a generator of  ${}_B\Mm^H$, and we can apply
\thref{8.6}.
\end{proof}


\begin{thebibliography}{99}
\bibitem{Bass}
H. Bass, Algebraic K-theory, Benjamin, New York, 1968.

\bibitem{Borceux}
F. Borceux, ``Handbook of categorical algebra I", {\sl Encyclopedia
Math. Appl.} {\bf 50}, Cambridge University Press, Cambridge, 1994.

\bibitem{Bourbaki}
N. Bourbaki, Alg\`ebre, Hermann, Paris, 1970.

\bibitem{Brzezinski02}
T. Brzezi\'nski, {The structure of corings. Induction functors,
Maschke-type theorem, and Frobenius and Galois properties},
{\sl Algebr. Representat. Theory} \textbf{5} (2002), 389--410.

\bibitem{BrzezinskiWisbauer}
T. Brzezi\'nski and R. Wisbauer, ``Corings and comodules",
{\sl London Math. Soc. Lect. Note Ser.} {\bf 309},
  Cambridge University Press, Cambridge, 2003.

\bibitem{Caenepeel03}
S. Caenepeel, {Galois corings from the descent theory point of view}, {\sl Fields
Inst. Comm.} {\bf 43} (2004), 163--186.

\bibitem{CMZ}
S. Caenepeel, G. Militaru and S. Zhu, ``Frobenius and separable
functors for generalized module categories and nonlinear
equations", {\sl Lect. Notes Math.} {\bf 1787}, Springer-Verlag,
Berlin, 2002.

\bibitem{CaenepeelR}
S. Caenepeel and \c S. Raianu, Induction functors for the Doi-Koppinen
unified Hopf modules, in ``Abelian groups and Modules", A. Facchini and
C. Menini (Eds.), Kluwer Academic Publishers, Dordrecht, 1995, p. 73--94.

\bibitem{DNR}
S. D\u asc\u alescu, C. N\u ast\u asescu and \c S.
Raianu, ``Hopf Algebras. An Introduction", {\sl Monographs Textbooks
Pure Appl. Math.} {\bf 235}, Marcel Dekker, New
York, 2001.

\bibitem{Doi}
Y. Doi, Unifying Hopf modules, {\sl J. Algebra} {\bf 153} (1992), 373--385.

\bibitem{DT}
Y. Doi and M. Takeuchi, {Hopf-Galois extensions of algebras, the Miyashita-Ulbrich
action, and Azumaya algebras}, {\sl J. Algebra} {\bf 121} (1989), 488--516.

\bibitem{KT}
Hopf algebras and Galois extensions of an algebra, {\sl Indiana Univ. Math. J.} {\bf 30}
(1981), 675--692.

\bibitem{Marcus}
A. Marcus, {Equivalences induced by graded bimodules},
{\sl Comm. Algebra} {\bf 26} (1998), 713--731.

\bibitem{Menini}
C. Menini and M. Zuccoli, {Equivalence theorems and Hopf-Galois extensions},  {\sl J. Algebra}  {\bf 194}  (1997),  245--274.

\bibitem{Mont}
S. Montgomery, ``Hopf algebras and their actions on rings", American
Mathematical
Society, Providence, 1993.

\bibitem{Schauenburg}
P. Schauenburg, Hopf bimodules over Hopf-Galois extensions, Miyashita-Ulbrich
actions, and monoidal center constructions, {\sl Comm. Algebra} {\bf 24}
(1996), 143-163.

\bibitem{Schneider0}
H.-J. Schneider, {Principal homogeneous
spaces for arbitrary Hopf algebras}, {\sl Israel J. Math.} {\bf 72}
(1990), 167--195.

\bibitem{Schneider1}
H.-J. Schneider, {Representation theory of Hopf Galois
extensions}, {\sl Israel J. Math.} {\bf 72} (1990), 196--231.

\bibitem{Schneider2}
H.-J. Schneider, \emph{Hopf Galois extensions,
crossed products, and Clifford theory}, Bergen, Jeffrey (ed.) et
al., Advances in Hopf algebras. Conference, August 10-14, 1992,
Chicago, IL, USA. New York, NY: Marcel Dekker. Lect. Notes Pure
Appl. Math. 158, 267-297 (1994).

\bibitem{Ulbrich}
K.-H. Ulbrich,  {On modules induced or coinduced from Hopf subalgebras}, {\sl Math. Scand.} {\bf 67}  (1990), 177--182.

\end{thebibliography}
\end{document}